\documentclass[11pt]{article}

\usepackage{amsmath,amssymb,amscd,amsthm}

\title{Some problems on mapping class groups and moduli space}

\author{Benson Farb \thanks{Supported in part by NSF grants DMS-9704640
and DMS-0244542.}}

\newtheorem{theorem}{Theorem}[section]
\newtheorem{proposition}[theorem]{Proposition}

\newtheorem{conjecture}[theorem]{Conjecture}
\newtheorem{question}[theorem]{Question}
\newtheorem{problem}[theorem]{Problem}

\newenvironment{example}
{\stepcounter{theorem}\bigskip\noindent{\bf Example
\arabic{section}.\arabic{theorem}.}}
{}

\def\bol{\bf}
\def\proof{{\bf {\medskip}{\noindent}Proof. }}

\def\endproof{$\diamond$ \bigskip}

\def\epar{\medskip}
\def\df{\em}
\def\title{\em}

\def\bar{\overline}

\newcommand\R{\mbox{\bf R}}
\newcommand\C{\mbox{\bf C}}
\newcommand\hyp{\mbox{\bf H}}
\newcommand\Z{\mbox{\bf Z}}
\newcommand\N{\mbox{\bf N}}
\newcommand\Q{\mbox{\bf Q}}

\newcommand\F{\mbox{\bf F}}

\DeclareMathOperator\kg{{\cal K}_g}
\DeclareMathOperator\Out{Out}
\DeclareMathOperator\Aut{Aut}
\DeclareMathOperator\rank{rank}
\DeclareMathOperator\spec{spec}
\DeclareMathOperator\cone{Cone}
\DeclareMathOperator\teich{Teich}
\DeclareMathOperator\Teich{\teich}
\DeclareMathOperator\MS{{\cal M}_g}
\DeclareMathOperator\Mod{Mod}
\DeclareMathOperator\BMod{BMod}
\DeclareMathOperator\Tor{{\cal I}}
\DeclareMathOperator\T{{\cal I}}
\DeclareMathOperator\K{{\cal K}}
\DeclareMathOperator\M{{\cal M}}
\DeclareMathOperator\U{{\cal U}}
\DeclareMathOperator\A{{\cal A}}
\DeclareMathOperator\Homeo{Homeo}
\DeclareMathOperator\Diff{Diff}
\DeclareMathOperator\BDiff{BDiff}

\DeclareMathOperator\Id{Id}

\DeclareMathOperator\Jac{Jac}
\DeclareMathOperator\Comm{Comm}
\DeclareMathOperator\Isom{Isom}

\DeclareMathOperator\SL{SL}
\DeclareMathOperator\PSL{PSL}
\DeclareMathOperator\SO{SO}
\DeclareMathOperator\SU{SU}
\DeclareMathOperator\GL{GL}
\DeclareMathOperator\Sp{Sp}

\DeclareMathOperator\e{ends}
\DeclareMathOperator\ent{ent}
\DeclareMathOperator\I{\Tor}
\DeclareMathOperator\IA{IA}
\DeclareMathOperator\Hom{Hom}
\DeclareMathOperator\Vol{Vol}
\DeclareMathOperator\cd{cd}
\DeclareMathOperator\image{image}

\newcommand{\margin}[1]{}

\begin{document}
\maketitle

\begin{abstract}
This paper presents a number of problems about mapping class groups and
moduli space.  The paper will appear in the book {\em Problems on
Mapping Class Groups and Related Topics}, ed. by B. Farb, {\em
Proc. Symp. Pure Math.} series, Amer. Math. Soc.
\end{abstract}

{\small \tableofcontents}

\section{Introduction}
This paper contains a biased and personal list of problems on mapping
class groups of surfaces.  One of the difficulties in this area has been
that there have not been so many easy problems.  One of my goals here is
to formulate a number of problems for which it seems that progress is
possible.  Another goal is the formulation of problems  which will force
us to penetrate more deeply into the structure of mapping class groups.
Useful topological tools have been developed, for example the Thurston
normal form, boundary theory, the reduction theory for subgroups, and
the geometry and topology of the complex of curves.  On the other hand
there are basic problems which seem beyond the reach of these methods.  One
of my goals here is to pose problems whose solutions might require new 
methods. 

\subsection{Universal properties of $\Mod_g$ and $\M_g$}

Let $\Sigma_g$ denote a closed, oriented surface of genus $g$, and let
$\Mod_g$ denote the group of homotopy classes of
orientation-preserving homeomorphisms of $\Sigma_g$.  
The mapping class group $\Mod_g$, along 
with its variations, derives much of its importance from its universal
properties.  Let me explain this for $g\geq 2$.  In this case, a
classical result of Earle-Eells gives that the identity component 
$\Diff^0(\Sigma_g)$ of $\Diff^+(\Sigma_g)$ is contractible.  Since 
$\Mod_g=\pi_0\Diff^+(\Sigma_g)$ by definition, we have a
homotopy equivalence of classifying spaces:
\begin{equation}
\label{eq:classifying1}
\BDiff^{+}(\Sigma_g)\simeq \BMod_g
\end{equation}

Let $\M_g$ denote the moduli space of Riemann surfaces.  
The group $\Mod_g$ acts properly discontinuously on the 
Teichm\"{u}ller space $\Teich_g$ of marked, genus $g$ Riemann surfaces.
Since $\Teich_g$ is contractible it 
follows that $\M_g$ is a $K(\Mod_g,1)$ space, i.e. it is homotopy
equivalent to the spaces in (\ref{eq:classifying1}).  From these
considerations it {\em morally} follows that, 
for any topological space $B$, we have the following bijections:

\begin{equation}
\label{eq:bijection1}
\left\{\begin{array}{c}
\mbox{Isomorphism classes}\\
\mbox{of $\Sigma_g$-bundles over}\\
\mbox{$B$}
\end{array}\right\}
\longleftrightarrow 
\left\{\begin{array}{c}
\mbox{Homotopy classes}\\
\mbox{of}\\
\mbox{maps $B\to\M_g$}
\end{array}\right\}
\longleftrightarrow 
\left\{\begin{array}{c}
\mbox{Conjugacy classes}\\
\mbox{of representations}\\
\mbox{$\rho:\pi_1B\to \Mod_g$}
\end{array}\right\}
\end{equation}

\bigskip

I use the term ``morally'' because (\ref{eq:bijection1}) is not exactly
true as stated.  For example, one can have two nonisomorphic 
$\Sigma_g$ bundles over $S^1$ with finite monodromy and with 
classifying maps $f:S^1\to
\Mod_g$ having the same image, namely a single point.  The problem here
comes from the torsion in 
$\Mod_g$.  This torsion prevents $\M_g$ from being a manifold; it is
instead an orbifold, and so we need to
work in the category of orbifolds.  This is a nontrivial issue which requires 
care.  There are two basic fixes to this problem.  First, one can 
simply replace $\M_g$ in (\ref{eq:bijection1}) 
with the classifying space $\BMod_g$.  Another option is to 
replace $\Mod_g$ in (\ref{eq:bijection1}) with any 
torsion-free subgroup $\Gamma<\Mod_g$ of finite index.  Then 
$\Gamma$ acts freely on $\Teich_g$ and the corresponding finite cover 
of $\M_g$ is a manifold.   In this case (\ref{eq:bijection1}) is true as
stated.  This torsion subtlety
will usually not have a major impact on our discussion, so we will 
for the most part ignore it.  This is fine on the level of rational
homology since the homotopy equivalences described above induce
isomorphisms: 

\begin{equation}
\label{eq:homol:isos}
H^\ast(\M_g,\Q)\approx H^\ast(\BDiff^{+}(\Sigma_g),\Q)\approx
H^\ast(\Mod_g,\Q)
\end{equation}

There is a unique complex orbifold structure on $\M_g$ with the property that 
these bijections carry over to the holomorphic
category.  This means that the manifolds are complex manifolds, the
bundles are non-isotrivial (i.e. the holomorphic structure of the fibers
is not locally constant, unless the map $B\to\M_g$ is trivial), 
and the maps are holomorphic.  For the third entry
of (\ref{eq:bijection1}), one must restrict to such conjugacy classes
with holomorphic representatives; many conjugacy classes do not have
such a representative.

For $g\geq 3$ 
there is a canonical $\Sigma_g$-bundle $\U_g$ over $\M_g$, called the
{\em universal curve} (terminology from algebraic geometry), for which
the (generic) fiber over any $X\in\M_g$ is the Riemann surface $X$.  The first
bijection of (\ref{eq:bijection1}) is realized concretely using $\U_g$.
For example given any smooth $f:B\to \M_g$, one simply pulls back the
bundle $\U_g$ over $f$ to give a $\Sigma_g$-bundle over $B$.  Thus
$\M_g$ plays the same role for surface bundles as the (infinite)
Grassmann manifolds play for vector bundles.  Again one needs to be
careful about torsion in $\Mod_g$ here, for example by passing to a
finite cover of $\M_g$.  For $g=2$ there are more serious problems.

An important consequence is the following.  Suppose one wants to
associate to every $\Sigma_g$-bundle a (say integral) 
cohomology class on the base of
that bundle, so that this association is {\em natural}, that is, 
it is preserved under pullbacks.  Then each such cohomology class must be the
pullback of some element of $H^\ast(\M_g,\Z)$.  In this sense the
classes in $H^\ast(\M_g,\Z)$ are universal.  After circle bundles, 
this is the next simplest nonlinear bundle theory.  
Unlike circle bundles, this study connects in a fundamental way
to algebraic geometry, among other things.

Understanding the sets in (\ref{eq:bijection1}) is interesting even in
the simplest cases.

\begin{example} ({\bf Surface bundles over $S^1$}). 
\label{example:sbundles}
Let $B=S^1$.  In this case (\ref{eq:bijection1}) states that the 
classfication of $\Sigma_g$-bundles over $S^1$, up to bundle
isomorphism, is equivalent to the classification of elements of $\Mod_g$ up to
conjugacy.  Now, a fixed $3$-manifold may fiber over $S^1$ in infinitely
many different ways, although there are finitely many fiberings with
fiber of fixed genus.  Since it is possible to compute these 
fiberings\footnote{This is essentially the fact that the Thurston norm is
computable.}, the homeomorphism problem for $3$-manifolds fibering over
$S^1$ can easily be reduced to solving the conjugacy problem for
$\Mod_g$.  This was first done by Hemion \cite{He}. 
\end{example}

\begin{example} ({\bf Arakelov-Parshin finiteness}). 
\label{example:geometric:shafarevich}
Now let $B=\Sigma_h$ for some $h\geq 1$, and consider the sets and 
bijections of (\ref{eq:bijection1}) in the holomorphic category.  
The Geometric Shafarevich Conjecture, proved by Arakelov and Parshin,
states that these sets are finite, and that the holomorphic 
map in each nontrivial homotopy class is unique.  As beautifully 
explained in \cite{Mc1}, from this
result one can derive (with a branched cover trick of Parshin) 
finiteness theorems for rational points on
algebraic varieties over function fields.
\end{example}

\medskip
\noindent
{\bf Remark on universality.} I would like to emphasize the following
point. While the existence of 
characteristic classes associated to {\em every} $\Sigma_g$-bundle 
is clearly the first case to look at, it seems that 
requiring such a broad form of 
universality is too constraining.  One reflection of this is the paucity
of cohomology of $\M_g$, as the Miller-Morita-Mumford conjecture (now
theorem due to Madsen-Weiss \cite{MW}) shows.  One problem is that the
requirement of naturality for {\em all} monodromies 
simply kills what are otherwise natural and common
classes.  Perhaps more natural would be to consider the characteristic
classes for $\Sigma_g$-bundles with torsion-free monodromy.  This would
lead one to understand the cohomology of various finite index subgroups
of $\Mod_g$.  

Another simple yet striking example of this phenomenon is Harer's
theorem that $$H^2(\M_g,\Q)=\Q$$ 
In particular the signature cocycle,
which assigns to every bundle $\Sigma_g\to M^4\to B$ the signature
$\sigma(M^4)$, is (up to a rational multiple) the only characteristic
class in dimension $2$.  When the monodromy representation
$\pi_1B\to\Mod_g$ lies in the (infinite index) 
Torelli subgroup $\T_g<\Mod_g$ (see
below), $\sigma$ is always zero, and so is useless.  However,  
there are infinitely many homotopy types of
surface bundles $M^4$ over surfaces with $\sigma(M^4)=0$; indeed such
families of examples can be taken to have monodromy in $\T_g$.  We 
note that there are no known elements of
$H^\ast(\Mod_g,\Q)$ which restrict to a nonzero element of
$H^\ast(\T_g,\Q)$. 

We can
then try to find, for example, characteristic classes for
$\Sigma_g$-bundles with monodromy lying in $\T_g$, and it is not hard to
prove that these are just pullbacks of classes in $H^\ast(\T_g,\Z)$.  In
dimensions one and two, for example, we obtain a large number of such
classes (see \cite{Jo1} and \cite{BFa2}, respectively).

I hope I have provided some motivation for understanding the cohomology
of subgroups of $\Mod_g$.  This topic is wide-open; we will discuss a
few aspects of it below.

\medskip
\noindent
{\bf Three general problems. }Understanding the theory of surface 
bundles thus leads to the following basic general problems.

\begin{enumerate}
\item For various finitely presented groups $\Gamma$, classify the 
representations $\rho:\Gamma\to\Mod_g$ up to conjugacy.

\item Try to find analytic and geometric structures on $\M_g$ in order
to apply complex and Riemannian geometry ideas to 
constrain possibilities on such $\rho$. 

\item Understand the cohomology algebra $H^\ast(\Gamma,K)$ for various
subgroups $\Gamma <\Mod_g$ and various modules $K$, and find topological and
geometric interpretations for its elements.  
\end{enumerate}

We discuss below some problems in these directions in
\S\ref{section:subgroup}, \S\ref{section:geomtop} and
\S\ref{section:torelli}, respectively.  One appealing 
aspect of such problems is that attempts at understanding
them quickly lead to ideas from combinatorial group theory, complex
and algebraic geometry, the theory of dynamical systems, low-dimensional
topology, symplectic representation theory, and more.  In addition to
these problems, we will be motivated here by the fact that $\Mod_g$ and
its subgroups provide a rich and important collection of examples to
study in combinatorial and geometric group theory; see
\S\ref{section:combinatorial:group:theory} below.

\medskip
\noindent
{\bf Remark on notational conventions. }
We will usually state conjectures and problems
and results for $\Mod_g$, that is, for closed surfaces.  This is simply
for convenience and simplicity; such conjectures and problems should
always be considered as posed for surfaces with boundary and punctures,
except perhaps for some sporadic, low-genus cases.  Similarly for other
subgroups such as the Torelli group $\T_g$.   Sometimes
the extension to these cases is straight-forward, but sometimes it
isn't, and new phenomena actually arise.  

\subsection{The Torelli group and associated subgroups}

One of the recurring objects in this paper will be the Torelli group.
We now briefly describe how it fits in to the general picture.

\medskip
\noindent
{\bf Torelli group. }Algebraic intersection number gives a symplectic form on 
$H_1(\Sigma_g,\Z)$.  This form is preserved by the natural action 
of $\Mod_g$.  The {\em Torelli group} $\T_g$ is defined to be the kernel
of this action.  We then have an exact sequence
\begin{equation}
\label{eq:torelli:main}
1\to \T_g\to \Mod_g\to \Sp(2g,\Z)\to 1
\end{equation}

The genus $g$ {\em Torelli space} is defined to be the quotient of
$\Teich_g$ by $\T_g$.  Like $\M_g$, this space has the appropriate 
universal mapping properties. However, the study of maps into Torelli
space is precisely {\em complementary} to the theory of holomorphic maps
into $\M_g$, as follows.  Any holomorphic map $f:B\to\M_g$ with
$f_\ast(B)\subseteq\T_g$, when composed with the (holomorphic) period mapping
$\M_g\to\A_g$ (see \S\ref{subsection:period:mapping} below), 
lifts to the universal cover $\widetilde{\A_g}$, which is
the Siegel upper half-space (i.e. the symmetric
space $\Sp(2g,\R)/\SU(g)$).  Since the domain is compact, the image of
this holomorphic lift is constant.  Hence $f$ is constant.

The study of $\T_g$ goes back to Nielsen (1919) and Magnus
(1936), although the next big breakthrough came in a series of
remarkable papers by Dennis Johnson in the late 1970's (see \cite{Jo1}
for a summary).  Still, many basic questions about $\T_g$ remain
open; we add to the list in \S\ref{section:torelli} below.

\medskip
\noindent
{\bf Group generated by twists about separating curves.} The 
{\em group generated by twists about separating curves},
denoted $\K_g$, is defined to be the subgroup $\Mod_g$ generated by the
(infinitely many) Dehn twists about separating (i.e. bounding) curves in
$\Sigma_g$.  The group $\K_g$ is sometimes called the 
{\em Johnson kernel} since Johnson proved that $\K_g$ is precisely the
kernel of the so-called Johnson homomorphism.  This group is a featured
player in the study of the Torelli group.  Its connection to
$3$-manifold theory begins with Morita's result that every integral
homology $3$-sphere comes from removing a handlebody component of some 
Heegaard embedding $h:\Sigma_g\hookrightarrow S^3$ and gluing it back 
to the boundary $\Sigma_g$ by an element of $\K_g$.  Morita then proves
the beautiful result that for a fixed such $h$, taking 
the Casson invariant of the resulting $3$-manifold actually gives a
{\em homomorphism} $\K_g\to\Z$.  This is a starting point for Morita's
analysis of the Casson invariant; see \cite{Mo1} for a summary.

\medskip
\noindent
{\bf Johnson filtration. }We now describe a filtration of $\Mod_g$ which
generalizes (\ref{eq:torelli:main}).  This filtration has become a basic
object of study.

For a group $\Gamma$ we inductively define $\Gamma_0:=\Gamma$ and
$\Gamma_{i+1}=[\Gamma,\Gamma_i]$.  The chain of subgroups
$\Gamma\supseteq \Gamma_1\supseteq \cdots $is the {\em lower
central series} of $\Gamma$. The group $\Gamma/\Gamma_i$ is $i$-step
nilpotent; indeed $\Gamma/\Gamma_i$ has the universal property that any
homomorphism from $\Gamma$ to any $i$-step nilpotent group factors
through $\Gamma/\Gamma_i$.  The sequence $\{\Gamma/\Gamma_i\}$ can be
thought of as a kind of Taylor series for $\Gamma$.  

Now let $\Gamma:=\pi_1\Sigma_g$.  It is a classical result of Magnus that 
$\bigcap_{i=1}^\infty \Gamma_i=1$, that is $\Gamma$ is {\em residually
nilpotent}.  Now $\Mod_g$ acts 
by outer automorphisms on $\Gamma$, and each of the subgroups $\Gamma_i$
is clearly characteristic.  We may thus define for each $k\geq 0$:
\begin{equation}
\label{eq:johnson:def}
\T_g(k):=\ker(\Mod_g\to \Out(\Gamma/\Gamma_k))
\end{equation}
That is, $\T_g(k)$ is just the subgroup of $\Mod_g$
acting trivially on the $k^{th}$ nilpotent quotient of
$\pi_1\Sigma_g$.  Clearly $\T_g(1)=\T_g$;  Johnson proved that
$\T_g(2)=\K_g$.  The sequence $\T_g=\T_g(1)\supset \T_g(2)\supset \cdots$
is called the {\em Johnson filtration}; it has also been called the {\em
relative weight filtration}.  This sequence forms a 
(but not {\em the}) lower central series for $\T_g$.  
The fact that $\bigcap_{i=1}^\infty \Gamma_i=1$ easily
implies that the $\Aut$ versions of the groups defined in 
(\ref{eq:johnson:def}) have trivial intersection.  The stronger fact that 
$\bigcap_{i=1}^\infty \T_g(i)=1$, so that $\T_g$ is
residually nilpotent, is also true, but needs some additional argument
(see \cite{BL}).

\bigskip
\noindent
{\bf Acknowledgements. }Thanks to everyone who sent in comments on
earlier drafts of this paper. Tara Brendle, Curt McMullen 
and Andy Putman made a number
of useful comments.  I am especially grateful to Chris Leininger, Dan
Margalit and Ben McReynolds, each of whom passed on to me extensive
lists of comments and corrections.  Their efforts have greatly improved
this paper.  Finally, I would like to express my deep appreciation
to John Stiefel for his constant friendship and support during the last
six months.

\section{Subgroups and submanifolds: Existence and
classification}
\label{section:subgroup}

It is a natural problem to classify subgroups $H$ of $\Mod_g$.  
By classify we mean to give an effective list of isomorphism or
commensurability types of subgroups, and also to determine all
embeddings of a given subgroup, up to conjugacy.  While 
there are some very general structure theorems,  one sees rather quickly
that the problem as stated is far too ambitous.  One thus begins with 
the problem of finding various invariants associated to subgroups, by
which we mean invariants of their isomorphism type, commensurability
type, or more extrinsic invariants which depend on the embedding $H\to
\Mod_g$.  One then tries to determine precisely which values of 
a particular invariant can occur, and
perhaps even classify those subgroups having a given value of the
invariant.  In this section we present a selection of such problems.

\medskip
\noindent
{\bf Remark on subvarieties. }The classification and existence problem for subgroups and
submanifolds of $\M_g$ can be viewed as algebraic and topological
analogues of the problem, studied by algebraic geometers, of
understanding (complete) subvarieties of $\M_g$.  There is an extensive
literature on this problem; see, e.g., \cite{Mor} for a survey.  To give
just one example of the type of problem studied, let 
$$c_g:=\max\{\dim_{\C}(V):
\mbox{$V$ is a complete subvariety of
$\M_g$}\}$$

The goal is to compute $c_g$.  This is a kind of measure of where $\M_g$
sits between being affine (in which case $c_g$ would be $0$) and
projective (in which case $c_g$ would equal $\dim_{\C}(\M_g)=3g-3$).  
While $\M_2$ is affine and so $c_2=0$, it is known for $g\geq 3$ that 
$1\leq c_g<g-1$; the lower bound is given by construction, the upper
bound is a well-known theorem of Diaz.

\subsection{Some invariants} 

The notion of {\em relative ends} $\e(\Gamma,H)$ 
provides a natural way to measure 
the ``codimension'' of a subgroup $H$ in a group $\Gamma$.  
To define $e(\Gamma,H)$, consider any proper, connected 
metric space $X$ on which $\Gamma$ acts properly and cocompactly by
isometries.  Then $\e(\Gamma, H)$ is defined to be the number of ends
of the quotient space $X/H$.

\begin{question}[Ends spectrum]
\label{question:ends}
What are the possibile values of $\e(\Mod_g,H)$ for finitely-generated
subgroups $H<\Mod_g$?
\end{question}

It is well-known that the moduli space $\M_g$ has one end.  The key
point of the proof is that the complex of curves is connected.
This proof actually gives more: any cover of $\M_g$ has one end; see,
e.g., \cite{FMa}.  However, I do not see how this fact directly gives
information about Question \ref{question:ends}.

\medskip
\noindent
{\bf Commensurators. }
Asking for two subgroups of a group to be conjugate is often too
restrictive a question.  A more robust notion is that of 
commensurability.  Subgroups $\Gamma_1,\Gamma_2$ of a group $H$ are {\em
commensurable} if there exists $h\in H$ such that $h\Gamma_1h^{-1}\cap
\Gamma_2$ has finite index in both $h\Gamma_1h^{-1}$ and in $\Gamma_2$.  
One then wants to classify subgroups up to commensurability; this
is the natural equivalence relation one studies in order to coarsify the
relation of ``conjugate'' to ignore finite index information.    
The
primary commensurability invariant for subgroups $\Gamma<H$ is 
the {\em commensurator of $\Gamma$ in $H$}, denoted $\Comm_H(\Gamma)$, 
defined as:
$$\Comm_H(\Gamma):=\{h\in H: h\Gamma h^{-1}\cap \Gamma \mbox{\ has finite
index in both $\Gamma$ and $h\Gamma h^{-1}$}\}$$
The commensurator has most commonly been studied for discrete subgroups
of Lie groups.  One of the most striking results about commensurators, 
due to Margulis, states that if $\Gamma$ is an irreducible lattice in a
semisimple\footnote{By {\em semisimple} we will always mean linear
semisimple with no compact factors.} Lie group $H$ then $[\Comm_H(\Gamma):\Gamma]=\infty$ if and
only if $\Gamma$ is arithmetic.  In other words, it is precisely the
arithmetic lattices that have infinitely many ``hidden symmetries''.

\begin{problem}
Compute $\Comm_{\Mod_g}(\Gamma)$ for various subgroups $\Gamma<\Mod_g$.  
\end{problem}

Paris-Rolfsen and Paris (see, e.g., \cite{Pa}) have proven that most subgroups
of $\Mod_g$ stabilizing a simple closed curve, or coming from the
mapping class group of a subsurface of $S$, are self-commensurating in
$\Mod_g$.  Self-commensurating subgroups, that is subgroups $\Gamma<H$
with $\Comm_H(\Gamma)=\Gamma$, are particularly important since the
finite-dimensional unitary dual of $\Gamma$ injects into the unitary
dual of $H$; in other words, any unitary representation of $H$ 
induced from a finite-dimensional irreducible unitary
representation of $\Gamma$ must itself be irreducible.

\medskip
\noindent
{\bf Volumes of representations. }Consider the general problem of
classifying, for a fixed finitely generated group $\Gamma$,
the set 
$${\cal X}_g(\Gamma):=\Hom(\Gamma,\Mod_g)/\Mod_g$$
of conjugacy classes of representations $\rho:\Gamma\to\Mod_g$.  Here
the representations $\rho_1$ and $\rho_2$ are conjugate if
$\rho_1=C_h\circ \rho_2$, where $C_h:\Mod_g\to\Mod_g$ is
conjugation by some $h\in\Mod_g$.  Suppose $\Gamma=\pi_1X$ where $X$
is, say, a smooth, closed $n$-manifold.  Since $\M_g$ is a classifying space 
for $\Mod_g$, we know that for each for each $[\rho]\in{\cal
X}_g(\Gamma)$ there exists a smooth map $f:X\to\M_g$ with
$f_\ast=\rho$, and that $f$ is unique up to homotopy.  

Each $n$-dimensional real cocycle $\xi$ on $\M_g$ then gives a
well-defined invariant 
$$\nu_\xi:{\cal X}_g(\Gamma)\to \R$$
defined by
$$\nu_\xi([\rho]):=\int_Xf^\ast \xi$$
It is clear that $\nu_\xi([\rho])$ does not depend on the choices,
and indeed depends only on the cohomology class of $\xi$.  As a case of
special interest, let $X$ be a $2k$-dimensional manifold and let 
$\omega_{\rm WP}$ denote the Weil-Petersson symplectic 
form  on $\M_g$.  Define the {\em complex $k$-volume} of
$\rho:\pi_1X\to \Mod_g$ to be 
$$\Vol_k([\rho]):=\int_Xf^\ast \omega_{\rm WP}^k$$

\begin{problem}[Volume spectrum]
\label{problem:volume:spec}
Determine for each $1\leq k\leq 3g-3$ the image of $\Vol_k:{\cal
X}_g(\Gamma)\to\R$.  Determine the union of all such images as $\Gamma$
ranges over all finitely presented groups.
\end{problem}
It would also be interesting to pose the same problem for
representations with special geometric constraints, for example those
with holomorphic or totally geodesic (with respect to a fixed metric)
representatives.  In particular, how do such geometric properties
constrain the set of possible volumes?  Note that Mirzakhani \cite{Mir}
has given recursive formulas for the Weil-Petersson volumes of moduli
spaces for surfaces with nonempty totally geodesic boundary.

\medskip
\noindent
{\bf Invariants from linear representations. }Each linear representation
$\psi:\Mod_g\to \GL_m(\C)$ provides us with many invariants for elements
of ${\cal X}_g(\Gamma)$, simply by composition with $\psi$ followed by
taking any fixed class function on $\GL_m(\C)$.  One can obtain
commensurability invariants for subgroups of $\Mod_g$ this way as well.
While no {\em faithful} $\psi$ is known for $g\geq 3$ (indeed the
existence of such a $\psi$ remains a major open problem) , there are
many such $\psi$ which give a great deal of information.  Some
computations using this idea can be found in \cite{Su2}.  I think
further computations would be worthwhile.

\subsection{Lattices in semisimple Lie groups}
\label{subsection:lattices}

While Ivanov proved that
$\Mod_g$ is not isomorphic to a lattice in a 
semisimple Lie group, a
recurring theme has been the comparison of algebraic properties of
$\Mod_g$ with such lattices and geometric/topological properties of
moduli space $\M_g$ with those of locally symmetric orbifolds.  The
question (probably due to Ivanov) then arose: ``Which lattices $\Gamma$
are subgroups of $\Mod_g$?''  This question arises from a slightly
different angle under the algebro-geometric guise of studying locally
symmetric subvarieties of moduli space; see \cite{Ha1}.  
The possibilities for such $\Gamma$ are highly constrained;
theorems of Kaimanovich-Masur,
Farb-Masur and Yeung (see \cite{FaM,Ye}), give the following.

\begin{theorem}
\label{theorem:rigid}
Let $\Gamma$ be an irreducible lattice in a semisimple Lie group $G\neq
\SO(m,1), \SU(n,1)$ with $m\geq 2, n\geq 1$.  Then any homomorphism
$\rho:\Gamma\to \Mod_g$ with $g\geq 1$ has finite image.
\end{theorem}

Theorem \ref{theorem:rigid} does not extend to the cases $G=\SO(m,1)$
and $G=\SU(n,1)$ in general since these groups admit lattices $\Gamma$
which surject to $\Z$.  Now let us restrict to the case of {\em injective}
$\rho$, so we are asking about which lattices $\Gamma$ occur as 
subgroups of $\Mod_g$.  As far as I know, here is what is currently
known about this question:
\begin{enumerate}
\item (Lattices $\Gamma<\SO(2,1)$): Each such $\Gamma$ has a finite
index subgroup which is free or $\pi_1\Sigma_h$ for some $h\geq 2$.  
These groups are plentiful in $\Mod_g$ and are discussed in more detail
in \S\ref{subsection:surfaces} below.

\item (Lattices $\Gamma<\SO(3,1)$): These exist
in $\Mod_{g,1}$ for $g\geq 2$ and in $\Mod_g$ for $g\geq 4$, 
by the following somewhat well-known construction.  
Consider the Birman exact sequence
\begin{equation}
1\to\pi_1\Sigma_g\to\Mod_{g,1}\stackrel{\pi}{\to}\Mod_g\to 1
\end{equation}

Let $\phi\in \Mod_g$ be a pseudo-Anosov homeomorphism, and let
$\Gamma_\phi<\Mod_{g,1}$ be the pullback under $\pi$ of the cyclic
subgroup of $\Mod_g$ generated by $\phi$.  The group
$\Gamma_\phi$ is isomorphic to the fundamental group of a $\Sigma_g$
bundle over $S^1$, namely the bundle obtained from $\Sigma_g\times
[0,1]$ by identifying $(x,0)$ with $(\phi(x),1)$.  
By a deep theorem of Thurston, such manifolds admit 
a hyperbolic metric, and so $\Gamma_\phi$ is a cocompact lattice in
$\SO(3,1)=\Isom^+(\hyp^3)$.  A branched covering trick of
Gonzalez-Diaz and Harvey (see \cite{GH}) can be used to find
$\Gamma_\phi$ as a subgroup of $\Mod_h$ for appropriate $h\geq 4$.  
A variation of the above can be used to find nonuniform lattices in 
$\SO(3,1)$ inside $\Mod_g$ for $g\geq 4$.

\item (Cocompact lattices $\Gamma<\SO(4,1)$): Recently John Crisp and I \cite{CF}
found one example of a cocompact lattice 
$\Gamma<\SO(4,1)=\Isom^+(\hyp^4)$ which embeds in $\Mod_g$ for all
sufficiently large $g$.  While we only know of one such $\Gamma$, it has
infinitely many
conjugacy classes in $\Mod_g$.  
The group $\Gamma$ is a right-angled Artin group, which is commensurable
with a group of reflections in the right-angled $120$-cell in $\hyp^4$. 

\item (Noncocompact lattices $\Gamma<\SU(n,1), n\geq 2$): These $\Gamma$
have nilpotent
subgroups which are not virtually abelian.  Since every nilpotent,
indeed solvable subgroup of $\Mod_g$ is virtually abelian, $\Gamma$ is
not isomorphic to any subgroup of $\Mod_g$. 
\end{enumerate}

Hence the problem of understanding which lattices in semisimple Lie
groups occur as subgroups of $\Mod_g$ comes down to the following.  

\begin{question}
\label{question:lattices2}
Does there exist some $\Mod_g, g\geq 2$ that contains a subgroup
$\Gamma$ isomorphic to a cocompact (resp. noncocompact) lattice in
$\SO(m,1)$ with $m\geq 5$ (resp. $m\geq 4$)? a
cocompact lattice in $\SU(n,1), n\geq 2$? Must there be only 
finitely many conjugacy classes of any such fixed $\Gamma$ in
$\Mod_g$?
\end{question}

In light of example (3) above, I would like to specifically ask: 
can $\Mod_g$ contain infinitely many isomorphism types of
cocompact lattices in $\SO(4,1)$?

Note that when $\Gamma$ is the fundamental group of a (complex)
algebraic variety $V$, then it is known that there can be at most finitely
many representations $\rho:\Gamma\to\Mod_g$ which have {\em holomorphic}
representatives, by which we mean the unique homotopy class of maps 
$f:V\to \M_g$ with $f_\ast=\rho$ contains a holomorphic map.  This
result follows from repeatedly taking hyperplane sections and finally
quoting the result for (complex) curves.  The result for these is 
a theorem of Arakelov-Parshin (cf. Example
\ref{example:geometric:shafarevich} above, and 
\S\ref{subsection:surfaces} below.)  

For representations which 
do not {\it a priori} have a holomorphic representative, one
might try to find a harmonic representative and then to prove a Siu-type
rigidity result to obtain a holomorphic representative.  One difficulty
here is that it is not easy to find harmonic representatives, since
(among other problems) every loop in $\M_g,g\geq 2$ can be freely
homotoped outside every compact set.  For recent progress, however, see 
\cite{DW} and the references contained therein.

\subsection{Surfaces in moduli space}
\label{subsection:surfaces}

Motivation for studying 
representations $\rho:\pi_1\Sigma_h\to \Mod_g$ for $g,h\geq
2$ comes from many directions.  These include: the analogy of
$\Mod_g$ with Kleinian groups (see, e.g., \cite{LR}); the fact that such
subgroups are the main source of locally symmetric families of 
Riemann surfaces (see \S\ref{subsection:lattices} above, and
\cite{Ha1}) ; and their appearance as a key piece of data 
in the topological classification of surface bundles over surfaces
(cf. (\ref{eq:bijection1}) above).  

Of course understanding such $\rho$ with holomorphic representatives is 
the Arakelov-Parshin Finiteness Theorem discussed in Example
\ref{example:geometric:shafarevich} above.  Holomorphicity is a key
feature of this result.  For example, one can prove (see, e.g.,
\cite{Mc1}) that there are finitely many
such $\rho$ with a holomorphic representative 
by finding a {\em Schwarz Lemma for $\M_g$}: any holomorphic
map from a compact hyperbolic surface into $\M_g$ endowed with the
Teichm\"{u}ller metric is distance decreasing.  The finiteness is also
just not true without the holomorphic assumption (see Theorem
\ref{tip:crisp} below).  We therefore want to recognize when a given
representation has a holomorphic representative.

\begin{problem}[Holomorphic representatives]
Find an algorithm or a group-theoretic invariant which determines or
detects whether or not a given representation
$\rho:\pi_1\Sigma_h\to \Mod_g$ has a holomorphic representative.
\end{problem}

Note that a necessary, but not
sufficient, condition for a representation
$\rho:\pi_1\Sigma_h\to\Mod_g$ to be holomorphic is that it be 
{\em irreducible}, i.e. there is no essential isotopy class of 
simple closed curve
$\alpha$ in $\Sigma_g$ such that
$\rho(\pi_1\Sigma_h)(\alpha)=\alpha$.  I believe it is not difficult
to give an algorithm to determine whether or not any given $\rho$ is
irreducible or not.

We would like to construct and classify (up to
conjugacy) such $\rho$.  We would also like to compute their associated
invariants, such as
$$\nu(\rho):=\int_{\Sigma_h} f^\ast\omega_{\rm WP}$$
where $\omega_{\rm WP}$ is the Weil-Petersson $2$-form on $\M_g$, and 
where $f:\Sigma_h\to\M_g$ is any map with $f_\ast=\rho$.  This would
give information on the signatures of surface bundles over surfaces, and
also on the Gromov co-norm of $[\omega_{\rm WP}]\in H^\ast(\M_g,\R)$.  

The classification question is basically impossible as stated, 
since e.g. surface
groups surject onto free groups, so it is natural to first restrict to 
injective $\rho$.  Using a technique of 
Crisp-Wiest, J. Crisp and I show in \cite{CF} 
that irreducible, injective $\rho$ are quite common.

\begin{theorem}
\label{tip:crisp}
For each $g\geq 4$ and each $h\geq 2$, there are 
infinitely many $Mod_g$-conjugacy classes of injective, 
irreducible representations $\rho:\pi_1\Sigma_h\to 
\Mod_g$.  One can take the images to lie inside 
the Torelli group $\T_g$.  Further, for any $n\geq 1$, one can take the 
images to lie inside the subgroup of $\Mod_g$ generated by $n^{\rm th}$ 
powers of all Dehn twists.  
\end{theorem}
                                                                            
One can try to use these representations, as well as those of \cite{LR},
to give new constructions of surface bundles over surfaces with small
genus base and fiber and nonzero signature.  The idea is that Meyer's 
signature cocycle is positively proportional to the Weil-Petersson $2$-form
$\omega_{WP}$, and one can actually explicitly integrate the pullback of
$\omega_{WP}$ under various representations, for example those glued
together using Teichm\"{u}ller curves.

Note that Theorem \ref{tip:crisp} provides
us with infinitely many topological types of surface bundles, each with 
irreducible, faithful monodromy, all having the same 
base and the same fiber, and all with signature zero.

\subsection{Normal subgroups}

It is a well-known open question to determine whether or not $\Mod_g$
contains a normal subgroup $H$ consisting of all pseudo-Anosov
homeomorphisms.  For genus $g=2$ Whittlesey \cite{Wh} found such an $H$;
it is an infinitely generated free group.  As far as I know this
problem is still wide open.

Actually, when starting to think about this problem I began to realize
that it is not easy to find {\em finitely generated} normal subgroups 
of $\Mod_g$ which are not commensurable with either $\T_g$ or $\Mod_g$.
There are
many normal subgroups of $\Mod_g$ which are not commensurable to either
$\T_g$ or to $\Mod_g$, most notably the terms $\T_g(k)$ of the Johnson
filtration for $k\geq 2$ and the terms 
of the lower central series of $\T_g$.  However, the former are 
infinitely generated and the latter are likely to be
infinitely generated; see Theorem
\ref{theorem:jf:notfg} and Conjecture \ref{conjecture:fglcc} below.  

\begin{question}[Normal subgroups]
\label{question:normal}
Let $\Gamma$ be a finitely generated normal subgroup of $\Mod_g$, where
$g\geq 3$.  Must $\Gamma$ be commensurable with $\Mod_g$ or with $\T_g$?
\end{question}

One way of constructing {\em infinitely generated} normal subgroups of
$\Mod_g$ is to take the group generated by the $n^{\rm th}$ powers of
all Dehn twists.  Another way is to take the normal closure $N_\phi$ of
a single element $\phi\in\Mod_g$.  It seems unclear how to determine the
algebraic structure of these $N_\phi$, in particular to determine
whether $N_\phi$ is finite index in $\Mod_g$, or in one of the
$\T_g(k)$.  The following is a basic test question.

\begin{question}
Is it true that, given any pseudo-Anosov $\phi\in\Mod_g$, there exists
$n=n(\phi)$ such that the normal closure of $\phi^n$ is free?
\end{question}

Gromov discovered the analogous phenomenon for elements of hyperbolic
type inside nonelementary word-hyperbolic groups; see \cite{Gro},
Theorem 5.3.E.  

One should compare Question \ref{question:normal} to the Margulis Normal
Subgroup Theorem (see, e.g. \cite{Ma}), which states that if $\Lambda$ is
any irreducible lattice in a real, linear semisimple Lie group 
with no compact factors and with $\R$-rank at 
least $2$, then any (not necessarily finitely generated) 
normal subgroup of $\Lambda$ is finite and central or has finite index in
$\Lambda$.  Indeed, we may apply this result to analyzing normal
subgroups $\Gamma$ of $\Mod_g$.  
For $g\geq 2$ the group $\Sp(2g,\Z)$ satisfies the hypotheses of
Margulis's theorem, and so the image $\pi(\Gamma)$ under
the natural representation $\pi:\Mod_g\to\Sp(2g,\Z)$ is normal, hence is
finite or finite index.  This proves the following.

\begin{proposition}[Maximality of Torelli]
\label{proposition:torelli:maximal}
Any normal subgroup $\Gamma$ of $\Mod_g$ containing $\T_g$ 
is commensurable either with $\Mod_g$ or with $\T_g$.
\end{proposition}

Proposition \ref{proposition:torelli:maximal} is a starting point for trying to
understand Question \ref{question:normal}.  Note too that 
Mess \cite{Me} proved that 
the group $\T_2$ is an infinitely generated free group, and so it
has no finitely generated normal subgroups.  Thus we know
that if $\Gamma$ is any finitely generated normal subgroup of $\Mod_2$,
then $\pi(\Gamma)$ has finite index in $\Sp(4,\Z)$.  This in turn gives
strong information about $\Gamma$; see \cite{Fa3}.

One can go further by considering the Malcev
Lie algebra ${\mathfrak t}_g$ of $\T_g$, computed by Hain in \cite{Ha3}; 
cf. \S\ref{section:liealg} below.  The normal 
subgroup $\Gamma\cap \T_g$ of $\T_g$ gives an $\Sp(2g,\Z)$-invariant 
subalgebra ${\mathfrak h}$ of ${\mathfrak t}_g$.  Let $H=H_1(\Sigma_g,\Z)$.  
The Johnson homomorphism
$\tau:\T_g\to \wedge^3H/H$ is equivariant with respect to the action of
$\Sp(2g,\Z)$, and $\wedge^3H/H$ is an irreducible $\Sp(2g,\Z)$-module.  It
follows that the first quotient in the lower central series of $\Gamma\cap
\T_g$ is either trivial or is all of $\wedge^3H/H$.  With more work, one
can extend the result of Proposition
\ref{proposition:torelli:maximal} from $\T_g$ to $\K_g=\T_g(2)$; 
see \cite{Fa3}.

\begin{theorem}[Normal subgroups containing $\K_g$]
Any normal subgroup $\Gamma$ of $\Mod_g$ containing $\K_g$ 
is commensurable to $\K_g$, $\T_g$ or $\Mod_g$.
\end{theorem}

One can continue this ``all or nothing image'' line of reasoning to
deeper levels of ${\mathfrak t}_g$.  Indeed, I believe one can completely 
reduce the classification of normal subgroups of $\Mod_g$, at least
those that contain some $\T_g(k)$, to some symplectic representation theory
problems, such as the following.

\begin{problem}
\label{problem:reps2}
For $g\geq 2$, determine the irreducible factors of the graded pieces of the 
Malcev Lie algebra ${\mathfrak t}_g$ of $\T_g$ as $\Sp$-modules.
\end{problem}

While Hain gives in \cite{Ha3} 
an explicit and reasonably simple presentation for 
${\mathfrak t}_g$ when $g>6$, Problem \ref{problem:reps2} still seems to be
an involved problem in classical representation theory.  

As Chris Leininger pointed out to me, all 
of the questions above have natural ``virtual versions''.  For
example, one can ask about the classification of 
normal subgroups of finite index
subgroups of $\Mod_g$.  Another variation is the classification of 
{\em virtually normal} subgroups of $\Mod_g$, that is, 
subgroups whose normalizers in $\Mod_g$ have finite index in $\Mod_g$.

\subsection{Numerology of finite subgroups}

The Nielsen Realization Theorem, due to Kerckhoff, states that any
finite subgroup $F<\Mod_g$ can be realized as a group of automorphisms
of some Riemann surface $X_F$. Here by {\em automorphism} group we mean
group of (orientation-preserving) isometries in some Riemannian metric,
or equivalently in the hyperbolic metric, or equivalently the group of
biholomorphic automorphisms of $X_F$.  An easy application of the
uniformization theorem gives these equivalences.  It is classical that
the automorphism group $\Aut(X)$ of a closed Riemann surface of genus
$g\geq 2$ is finite.  Thus the study of finite subgroups of $\Mod_g,
g\geq 2$ reduces to the study of automorphism groups of Riemann
surfaces.

Let $N(g)$ denote the largest
possible order of $\Aut(X)$ as $X$ ranges over all genus $g$ surfaces.  
Then for $g\geq 2$ it is known that
\begin{equation}
\label{eq:auto:bounds}
8(g+1)\leq N(g)\leq 84(g-1)
\end{equation}
the lower bound due to Accola and Maclachlan (see, e.g., \cite{Ac}); 
the upper due to Hurwitz.
It is also known that each of these bounds is achieved for infinitely
many $g$ and is not achieved for infinitely many $g$.  
There is an extensive literature seeking to compute $N(g)$ and
variations of it.  The most significant achievement in this direction is
the following result of M. Larsen \cite{La}.

\begin{theorem}[Larsen]
\label{theorem:larsen}
Let $H$ denote the set of integers $g\geq 2$ such that there exists at
least one compact Riemann surface $X_g$ of genus $g$ with
$|\Aut(X_g)|=84(g-1)$.  Then the series $\sum_{g\in H}g^{-s}$ converges
absolutely for $\Re(s)>1/3$ and has a singularity at $s=1/3$.
\end{theorem}

In particular, the $g$ for which $N(g)=84(g-1)$ 
occur with the same frequency as
perfect cubes.  It follows easily from the Riemann-Hurwitz formula that
the bound $84(g-1)$ is achieved precisely for those surfaces which
isometrically (orbifold) cover the $(2,3,7)$ orbifold.  Thus the problem of
determining which genera realize the $84(g-1)$ bound comes down to 
figuring out the possible (finite) indices of subgroups, or what is 
close to the same thing finite quotients, of the $(2,3,7)$ orbifold
group. Larsen's argument uses in a fundamental way 
the classification of finite simple groups.

\begin{problem}
Give a proof of Theorem \ref{theorem:larsen} which does not depend on
the classification of finite simple groups.
\end{problem}

To complete the picture, one would like to understand the frequency of
those $g$ for which the lower bound in (\ref{eq:auto:bounds}) 
occurs. 

\begin{problem}[Frequency of low symmetry]
\label{problem:accola}
Let $H$ denote the set of integers $g\geq 2$ such that $N(g)=8(g+1)$.  
Find the $s_0$ for which the 
series $\sum_{g\in H}g^{-s}$ converges
absolutely for the real part $\Re(s)$ of $s$ satisfying 
$\Re(s)>s_0$, and has a singularity at $s=s_0$.
\end{problem}

There are various refinements and variations on Problem
\ref{problem:accola}.  For example, Accola proves in \cite{Ac} that when
$g$ is divisible by $3$, then $N(g)\geq 8(g+3)$, with the bound attained
infinitely often.  One can try to build on this for other $g$, and can
also ask for the frequency of this occurence.

One can begin to refine Hurwitz's Theorem by asking for bounds of orders
of groups of automorphisms which in addition satisfy various algebraic
constraints, such as being nilpotent, being 
solvable, being a $p$-group, etc.  There already exist a number of
theorems of this sort.  For example, Zomorrodian \cite{Zo} proved
that if $\Aut(X_g)$ is nilpotent then it has order at most $16(g-1)$,
and if this bound is attained then $g-1$ must be a power of
$2$. One can also ask for lower bounds in this context.  As these kinds
of bounds are typically attained and not attained for infinitely many
$g$, one then wants to solve the following.

\begin{problem}[Automorphism groups with special properties]
Let $P$ be a property of finite groups, for example being nilpotent,
solvable, or a $p$-group.  Prove a version of Larsen's theorem which
counts those $g$ for which the upper bound of $|\Aut(X_g)|$ is realized
for some $X_g$ with $Aut(X_g)$ having $P$. Similarly for lower bounds.
Determine the least $g$ for which each given bound is realized.
\end{problem}

Many of the surfaces realizing the extremal bounds in all of the above
questions are {\em arithmetic}, that is they are quotients of $\hyp^2$
by an arithmetic lattice.  Such lattices are well-known to have special
properties, in particular they have a lot of symmetry.  On the other
hand arithmetic surfaces are not typical.  Thus to understand the
``typical'' surface with symmetry, the natural problem is the following.

\begin{problem}[Nonarithmetic extremal surfaces]
\label{problem:nonarithmetic}
Give answers to all of the above problems on automorphisms of Riemann
surfaces for the collection of {\em non-arithmetic} surfaces.  
For example find bounds on orders of automorphism groups which are
nilpotent, solvable, $p$-groups, etc.  Prove that these bounds are
sharp for infinitely many $g$.  Determine 
the frequency of those $g$ for which such bounds are sharp.  Determine
the least genus for which the bounds are sharp.
\end{problem}

The model result for these kind of problems is the following.

\begin{theorem}[Belolipetsky \cite{Be}]
\label{theorem:belo}
Let $X_g$ be any non-arithmetic Riemann surface of genus $g\geq 2$.
Then 
$$|\Aut(X_g)|\leq \frac{156}{7}(g-1)
$$
Further, this bound is sharp for infinitely many $g$; the least such is
$g=50$.  
\end{theorem}

The key idea in the proof of Theorem \ref{theorem:belo} is the
following.  One considers the quotient 
orbifold $X_g/\Aut(X_g)$, and computes via the Riemann-Hurwitz formula
that it is the quotient of $\hyp^2$ by a triangle group.  Lower bounds 
on the area of this orbifold give upper bounds on $|\Aut(X_g)|$; for
example the universal lower bound of $\pi/21$ for the area of every 
$2$-dimensional 
hyperbolic orbifold gives the classical Hurwitz $84(g-1)$ theorem.  
Now Takeuchi classified all arithmetic triangle groups; they are given by
a finite list.  One can then use this list to refine the usual
calculations with the Riemann-Hurwitz formula to give results such as 
Theorem \ref{theorem:belo}; see \cite{Be}.  This idea should also be 
applicable to Problem \ref{problem:nonarithmetic}.  

I would like to mention a second instance of the theme of playing off
algebraic properties of automorphism groups versus the numerology of
their orders.  One can prove that the maximal possible order of an
automorphism of a Riemann surface $X_g$ is $4g+2$.  This bound is easily
seen to be achieved for each $g\geq 2$ by considering the rotation of the
appropriately regular hyperbolic $(4g+2)$-gon.  
Kulkarni \cite{Ku} proved that there is a {\em
unique} Riemann surface $W_g$ admitting such an automorphism.  Further, he
proved that $\Aut(W_g)$ is cyclic, and he gave the equation describing
$W_g$ as an algebraic curve.

\begin{problem}[Canonical basepoints for $\M_g$]
Find other properties of automorphisms or automorphism groups that
determine a unique point of $\M_g$.  For example, is there a unique
Riemann surface of genus $g\geq 2$ whose automorphism group is
nilpotent, and is the largest possible order among nilpotent
automorphism groups of genus $g$ surfaces?
\end{problem}

A {\em Hurwitz surface} is a hyperbolic surface attaining
the bound $84(g-1)$.  As the quotient of such a surface by its 
automorphism group is the $(2,3,7)$ orbifold, which has a unique
hyperbolic metric, it follows that for each $g\geq 2$ there are finitely
many Hurwitz surfaces.  As these are the surfaces of maximal symmetry, 
it is natural to ask precisely how many there are.

\begin{question}[Number of Hurwitz surfaces]
Give a formula for the number of Hurwitz surfaces of genus
$g$. What is the frequency of those $g$ for which there is a unique
Hurwitz surface?
\end{question}

\section{Combinatorial and geometric group theory of $\Mod_g$}
\label{section:combinatorial:group:theory}

Ever since Dehn, $\Mod_g$ has been a central example in combinatorial
and geometric group theory.  One reason for this is that $\Mod_g$ lies
at a gateway:  on one side are matrix groups and groups naturally
equipped with a geometric structure (e.g. hyperbolic geometry); on the
other side are groups given purely combinatorially.  In this section we 
pose some problems in this direction.

\subsection{Decision problems and almost convexity}

\medskip
\noindent
{\bf Word and conjugacy problems. }
Recall that the {\em word problem} for a finitely presented group
$\Gamma$ asks for an algorithm which takes as input any word $w$ in a fixed
generating set for $\Gamma$, and as output tells whether or not $w$
is trivial.  The {\em conjugacy problem} for $\Gamma$ asks for an
algorithm which takes as input two
words, and as output tells whether or not these words represent
conjugate elements of $\Gamma$.  

There is some history to the word and conjugacy problems for $\Mod_g$,  
beginning with braid groups.  These problems have topological
applications; for example the conjugacy problem for $\Mod_g$ is one of
the two main ingredients one needs to solve the homeomorphism problem
for $3$-manifolds fibering over the circle\footnote{The second
ingredient, crucial but ignored by some authors, is the computability of
the Thurston norm.}.  Lee Mosher proved in \cite{Mos} that $\Mod_g$ is
automatic.  From this it follows that there is an $O(n^2)$-time
algorithm to solve the word problem for $\Mod_g$; indeed there is an
$O(n^2)$-time algorithm which puts each word in a fixed generating set
into a unique normal form. However, the following is still open.

\begin{question}[Fast word problem]
Is there a sub-quadratic time algorithm to solve the word problem in
$\Mod_g$?
\end{question}

One might guess that $n\log n$ is possible here, as there is such an
algorithm for certain relatively (strongly) hyperbolic groups 
(see \cite{Fa3}), and mapping class groups are at least weakly
hyperbolic, as proven by Masur-Minsky (Theorem 1.3 of \cite{MM1}).

The conjugacy problem for $\Mod_g$ is harder.  The original algorithm 
of Hemion \cite{He} seems to give no reasonable (even exponential) time
bound. One refinement of the problem would be to prove that 
$\Mod_g$ is {\em biautomatic}.  
However, even a biautomatic structure gives
only an exponential time algorithm to solve the conjugacy problem.
Another approach to solving the conjugacy problem is the following. 

\begin{problem}[Conjugator length bounds]
\label{problem:conj:length}
Prove that there 
exist constants $C,K$, depending only on $S$, so that if $u,v\in
\Mod_g$ are conjugate, then there exists $g\in \Mod_g$ with $||g||\leq
K\max\{||u||,||v||\}+C$ so that $u=gvg^{-1}$.
\end{problem}

Masur-Minsky (\cite{MM2}, Theorem 7.2) 
solved this problem in the case where $u$ and
$v$ are pseudo-Anosov; their method of hierarchies seems quite
applicable to solving Problem \ref{problem:conj:length} in the general
case.  While interesting in its own right, even the solution to Problem
\ref{problem:conj:length} would not answer the following basic problem.

\begin{problem}[Fast conjugacy problem]
\label{problem:conj:fast}
Find a polynomial time algorithm to solve the conjugacy problem in
$\Mod_g$. Is there a quadratic time algorithm, as for the word problem?
\end{problem}

As explained in the example on page \pageref{example:sbundles}, a 
solution to Problem \ref{problem:conj:fast} would be a major step 
in finding a polynomial time algorithm to solve the homeomorphism
problem for $3$-manifolds that fiber over the circle.

\medskip
\noindent
{\bf Almost convexity. }
In \cite{Ca} Cannon initiated the beautiful theory of almost convex
groups.  A group $\Gamma$ with generating set $S$ is {\em almost convex}
if there exists $C>0$ so that for each $r>0$, and for 
any two points $x,y\in \Gamma$ on the sphere of radius $r$ in $\Gamma$
with $d(x,y)=2$, there exists a path $\gamma$ of length at most $C$ 
connecting $x$ to $y$ and 
lying completely inside the ball of radius $r$ in $\Gamma$.  
There is an obvious generalization
of this concept from groups to spaces.  One strong consequence of the almost
convexity of $\Gamma$ 
is that for such groups one can recursively build the Cayley graph of
$\Gamma$ near any point $x$ in the $n$-sphere of $\Gamma$ by only 
doing a local computation involving elements of $\Gamma$ lying
(universally) close to $x$; see \cite{Ca}.  In particular one can
build each $n$-ball, $n\geq 0$, and so solve the word problem in an
efficient way. 

\begin{question}[Almost convexity]
\label{question:ac1}
Does there exist a finite generating set for $\Mod_g$ for which it is 
almost convex?
\end{question}

One would also like to know the answer to this question for various
subgroups of $\Mod_g$.  
Here is a related, but different, basic question about the geometry of
Teichm\"{u}ller space.

\begin{question}
\label{question:ac2}
Is $\Teich(\Sigma_g)$, endowed with the Teichm\"{u}ller metric, 
almost convex?  
\end{question}

Note that Cannon proves in \cite{Ca} that fundamental groups of
closed, negatively curved manifolds are almost convex with respect to
any generating set.  He conjectures that this should generalize both 
to the finite volume case and to the nonpositively curved case.  

\subsection{The generalized word problem and distortion}
\label{subsection:distortion}

In this subsection we pose some problems relating to the ways in which 
subgroups embed in $\Mod_g$.

\medskip
\noindent
{\bf Generalized word problem. }
Let $\Gamma$ be a finitely presented group with finite generating set 
$S$, and 
let $H$ be a finitely generated subgroup of $\Gamma$.  
The {\df generalized word problem}, or GWP, for $H$ in $\Gamma$,  
asks for an algorithm which takes as input an element of 
the free group on $S$, and as output tells 
whether or not this element represents an element of $H$.  
When $H$ is the trivial subgroup then this is simply the word problem
for $\Gamma$.  The group $\Gamma$ is said to have 
{\df solvable generalized word problem} if 
the GWP is solvable for every finitely generated subgroup $H$ in $\Gamma$.
\epar

The generalized word problem, also called the {\df membership, occurrence} 
or {\df Magnus problem}, was formulated by K. Mihailova \cite{Mi} in 
1958, but special cases had already been studied by 
Nielsen \cite{Ni} and Magnus \cite{Ma}.  Mihailova \cite{Mi} found a
finitely generated subgroup of a product $F_m\times F_m$ of free groups
which has unsolvable generalized word problem.

Now when $g\geq 2$, the group $\Mod_g$ contains a product of
free groups $F_m\times F_m$ (for any $m>0$).  Further, one can find
copies of $F_m\times F_m$ such that the generalized word problem for
these subgroups is solvable inside $\Mod_g$.  To give one concrete
example, simply divide $\Sigma_g$ into two subsurfaces $S_1$ and $S_2$
which intersect in a (possibly empty) collection of curves.  Then take,
for each $i=1,2$, a pair $f_i,g_i$ of independent pseudo-Anosov
homeomorphisms.  After perhaps taking powers of $f_i$ and $g_i$ if
necessary, the group generated by $\{f_1,g_1,f_2,g_2\}$ will be the
group we require; one can pass to finite index subgroups if one wants
$m>2$.  The point here is that such subgroups are not
distorted (see below), as they are convex cocompact in the subgroups
$\Mod(S_i)$ of $\Mod_g$, and these in turn are not distorted in
$\Mod_g$.  It follows that the generalized word problem for Mihailova's 
subgroup $G<F_m\times F_m$ is not solvable in
$\Mod_g$.

\begin{problem}[Generalized word problem]
\label{problem:gwp1}
Determine the subgroups $H$ in $\Mod_g$ for which the generalized word
problem is solvable.  Give efficient algorithms to solve the generalized
word problem for these subgroups.  Find the 
optimal time bounds for such algorithms.    
\end{problem}

Of course this problem is too broad to solve in complete generality;
results even in special cases might be interesting.  Some solutions to
Problem \ref{problem:gwp1} are given by Leininger-McReynolds in
\cite{LM}.  We also note that Bridson-Miller (personal communication), 
extending an old result of
Baumslag-Roseblade, have proven that any product of finite rank free
groups has solvable generalized word problem with respect to 
any {\em finitely presented} subgroup.  In light of this result, 
it would be interesting to determine whether or not there is a finitely
presented subgroup of $\Mod_g$ with respect to which the generalized
word problem is not solvable.

\medskip
\noindent
{\bf Distortion and quasiconvexity. }
There is a refinement of Problem \ref{problem:gwp1}.  Let $H$ be a
finitely generated subgroup of a finitely generated group $\Gamma$.  
Fix finite generating sets on both $H$ and $\Gamma$.  This choice gives a 
{\em word metric} on both $H$ and $\Gamma$, 
where $d_\Gamma(g,h)$ is defined to be the minimal number of
generators of $\Gamma$ needed to represent $gh^{-1}$. 
Let $\bol N$ denote the natural numbers.  
We say that a function $f:{\bol N}\longrightarrow {\bol N}$ is a 
{\df distortion function} for $H$ in $\Gamma$ if for every word $w$ in
the generators of $\Gamma$, if $w$ represents an element $\bar{w}\in H$ then 
$$d_H(1,\bar{w})\leq f(d_\Gamma(1,\bar{w}))$$

In this case we also say that ``$H$ has distortion 
$f(n)$ in $\Gamma$.''  It is easy to see that the growth type of $f$,
i.e. polynomial, exponential, etc., does not depend on the choice of
generators for either $H$ or $\Gamma$.  It is also easy to see that 
$f(n)$ is constant if and only if $H$ is 
finite; otherwise $f$ is at least linear.  It is proved in \cite{Fa1}
that, for a group $\Gamma$ with solvable word problem, 
the distortion of $H$ in $\Gamma$ is recursive if and only if $H$ has
solvable generalized word problem in $\Gamma$.  For some concrete examples, 
we note that the center of
the $3$-dimensional integral Heisenberg group has quadratic distortion;
and the cyclic group generated by $b$ in the group 
$\langle a,b: aba^{-1}=b^2\rangle$ has
exponential distortion since $a^nba^{-n}=b^{2^n}$.

\begin{problem}[Distortion]
Find the possible distortions of subgroups in $\Mod_g$.  In particular,
compute the distortions of $\T_g$.  Determine the
asymptotics of the distortion of $\T_g(k)$ as $k\to\infty$.  
Is there a subgroup $H<\Mod_g$ that
has precisely polynomial distortion of degree $d>1$?
\end{problem}

There are some known results on distortion of
subgroups in $\Mod_g$.  Convex cocompact subgroups (in the sense of
\cite{FMo}) have linear distortion in $\Mod_g$; there are many such
examples where $H$ is a free group.  Abelian subgroups of $\Mod_g$ 
have linear distortion 
(see \cite{FLMi}), as do subgroups corresponding to mapping
class groups of subsurfaces (see \cite{MM2,Ham}).

I would guess that $\I_g$ has exponential distortion
in $\Mod_g$.  A first step to the question of how 
the distortion of $\T_g(k)$ in $\Mod_g$
behaves as $k\to
\infty$ would be to determine the distortion
of $\T_g(k+1)$ in $\T_g(k)$.  The ``higher Johnson
homomorphisms'' (see, e.g., \cite{Mo3}) might be useful here.

A stronger notion than linear distortion is that of quasiconvexity. 
Let $S=S^{-1}$ be a fixed generating set for a group $\Gamma$, and let 
$\pi:S^\ast\to \Gamma$ be the natural surjective homomorphism 
from the free monoid on $S$ to $\Gamma$ sending a word to the group
element it represents.  Let $\sigma:\Gamma\to S^\ast$ be a (perhaps
multi-valued) section of $\pi$; that is, $\sigma$ is just a choice of
paths in $\Gamma$ from the origin to each $g\in\Gamma$.  We say that a
subgroup $H<\Gamma$ is {\em quasiconvex} (with respect to $\sigma$) 
if there exists $K>0$ so that for each $h\in H$, each path 
$\sigma(h)$ lies in the $K$-neighborhood of
$H$ in $\Gamma$. Quasiconvexity is a well-known and basic notion in
geometric group theory.  It is easy to see  
that if $H$ is quasiconvex with respect to
some collection of quasigeodesics then $H$ has linear distortion in 
$\Mod_g$ (see \cite{Fa1}). 

\begin{problem}
Determine which subgroups of $\Mod_g$ are quasiconvex with respect to
some collection of geodesics.
\end{problem}

This question is closely related to, but different than, the question 
of {\em convex cocompactness} of subgroups of $\Mod_g$, as defined in 
\cite{FMo}, since the embedding of $\Mod_g$ in $\Teich_g$ via any orbit
is exponentially distorted, by \cite{FLMi}.

\subsection{Decision problems for subgroups}

As a collection of groups, how rich and varied 
can the set of subgroups of $\Mod_g$ be?  One instance of this general 
question is the following. 

\begin{question}
Does every finitely presented subgroup $H<\Mod_g$ have solvable
conjugacy problem? is it combable? 
automatic?
\end{question}

Note that every finitely-generated subgroup of a group with solvable
word problem has solvable word problem.  The same is not true for the
conjugacy problem: there are subgroups of $\GL(n,\Z)$ 
with unsolvable conjugacy problem; see 
\cite{Mi}.  

It is not hard to see that $\Mod_g$, like $\GL(n,\Z)$, has finitely many
conjugacy classes of finite subgroups.  However, we pose the
following.

\begin{problem}
\label{problem:bridson:conj}
Find a finitely presented subgroup $H<\Mod_g$ for which there are
infinitely many conjugacy classes of finite subgroups in $H$.  
\end{problem}

The motivation for this problem comes from a corresponding example, due
to Bridson \cite{Br}, of such an $H$ in $\GL(n,\Z)$.  One might to 
solve Problem \ref{problem:bridson:conj} by 
extending Bridson's construction to $\Sp(2g,\Z)$, pulling back such
an $H$, and also noting that the natural map $\Mod_g\to\Sp(2g,\Z)$ is
injective on torsion.   

Another determination of the variety of subgroups of $\Mod_g$ is the
following.  Recall that the {\em isomorphism problem} for a collection
${\cal S}$ of finitely presented groups asks for an algorithm which
takes as input two presentations for two elements of ${\cal S}$ and as
output tells whether or not those groups are isomorphic.  

\begin{question}[Isomorphism problem for subgroups]
Is the isomorphism problem for the collection of finitely presented
subgroups of $\Mod_g$ solvable?  
\end{question}

Note that the isomorphism problem is not solvable
for the collection of all finitely generated linear groups, nor is it
even solvable for the collection of finitely generated subgroups of
$\GL(n,\Z)$.  

There are many other algorithmic questions one can ask; 
we mention just one more.

\begin{question}
Is there an algorithm to decide whether or not a given subgroup
$H<\Mod_g$ is freely indecomposable? Whether or not $H$ splits 
over $\Z$?
\end{question}

\subsection{Growth and counting questions}

Recall that the {\em growth series} of a group $\Gamma$ with respect to
a finite generating set $S$ is defined to be the power series  
\begin{equation}
\label{eq:growth:series}
f(z)=\sum_{i=0}^\infty c_iz^i
\end{equation}
where $c_i$ denotes the cardinality $\#B_\Gamma(i)$ 
of the ball of radius $i$ in
$\Gamma$ with respect to the word metric induced by $S$. We say that
$\Gamma$ has {\em rational growth} (with respect to $S$) if $f$ is a
rational function, that is the quotient of two polynomials.  This is
equivalent to the existence of a linear recurrence relation among the
$c_i$; that is, there exist $m>0$ real 
numbers $a_1,\ldots ,a_m\geq 0$ so that for each $r$:
$$c_r=a_1c_{r-1}+\cdots +a_mc_{r-m}$$

Many groups have rational growth with respect to various (sometimes
every) generating sets.  Examples include word-hyperbolic groups,
abelian groups, and Coxeter groups.  See, e.g., \cite{Harp2} 
for an introduction to the theory of growth series.  

In \cite{Mos}, Mosher constructed an automatic structure for $\Mod_g$.
This result suggests that the following might have a positive answer.

\begin{question}[Rational growth]
\label{question:rational:growth}
Does $\Mod_g$ have rational growth function with respect to some set 
of generators?  with respect to every set of generators?  
\end{question}

Of course one can also ask the same question for any 
finitely generated subgroup of $\Mod_g$, for example $\T_g$.  
Note that the existence of
an automatic structure is not known to imply rationality of growth (even
for one generating set); one needs in addition the property that the
automatic structure consist of geodesics.  Unfortunately Mosher's
automatic structure does not satisfy this stronger condition.  

It is natural to ask for other recursive patterns in the Cayley graph of
$\Mod_g$.  To be more precise, let P denote a property that elements of
$\Mod_g$ might or might not have.  For example, P might be the property 
of being finite order, of lying in a fixed subgroup $H<\Mod_g$, or of
being pseudo-Anosov.  Now let
$$c_P(r)=\#\{B_{\Mod_g}(r)\cap\{x\in\Mod_g:\mbox{$x$ has P}\}\}$$

We now define the {\em growth series for the property P}, with respect
to a fixed generating set for $\Gamma$, to be the power series
$$
f_P(z)=\sum_{i=0}^\infty c_P(i)z^i
$$

\begin{question}[Rational growth for properties]
For which properties P is the function $f_P$ is rational?  
\end{question}

\medskip
\noindent
{\bf Densities. }For any subset $S\subset\Mod_g$, it is natural 
to ask how common elements of $S$ are in $\Mod_g$.  There are
various ways to interpret this question, and the answer likely depends
in a strong way on the choice of interpretation\footnote{For a
wonderful discussion of this kind of 
issue, see Barry Mazur's article \cite{Maz}.}.   
One way to formalize this is via the {\em density}
$d(S)$ of $S$ in $\Mod_g$, where
\begin{equation}
\label{eq:density}
d(S)=\lim_{r\to\infty}\frac{\#[ B(r)\cap S]}{\# B(r)}
\end{equation}
where $B(r)$ is the number of elements of $\Mod_g$ in the ball of radius
$r$, with respect to a fixed set of generators.  While for subgroups
$H<\Mod_g$ the number $d(H)$ itself may
depend on the choice of generating sets for $H$ and $\Mod_g$, it is not
hard to see that the (non)positivity of $d(H)$ does not depend on the
choices of generating sets. 

As the denominator and (typically) the
numerator in (\ref{eq:density}) are exponential, one expects that
$d(S)=0$ for most $S$.  Thus is it natural to replace both the numerator
and denominator of (\ref{eq:density}) with their logarithms; we denote the
corresponding limit as in (\ref{eq:density}) by $d_{\log}(S)$, and we
call this the {\em logarithmic density} of $S$ in $\Mod_g$.

The following is one interpretation of a folklore conjecture.

\begin{conjecture}[Density of pseudo-Anosovs]
\label{conjecture:pseudos:dense}
Let ${\cal P}$ denote the set of pseudo-Anosov elements of $\Mod_g$.  
Then $d({\cal P})=1$.
\end{conjecture}

J. Maher \cite{Mah} 
has recently proven that a random walk on $\Mod_g$ lands on a
pseudo-Anosov element with probability tending to one as the length of
the walk tends to infinity.  I. Rivin \cite{Ri} 
has proven that a random (in a certain specific sense) element of
$\Mod_g$ is pseudo-Anosov by proving a corresponding result for
$\Sp(2g,\Z)$.  While the methods in \cite{Mah} and \cite{Ri} may be
relevant, Conjecture \ref{conjecture:pseudos:dense} 
does not seem to follow directly from these results.   As another test 
we pose the following.

\begin{conjecture}
\label{conjecture:dtor}
$d(\Tor_g)=0$.
\end{conjecture}

Even better would be to determine $d_{\log}(\T_g)$.  
Conjecture \ref{conjecture:dtor} would
imply that $d(\T_g(m))=0$ for each $m\geq 2$. It is not hard to see that $\T_g(m)$ has exponential growth for
each $g\geq 2, m\geq 1$, and one wants to understand
how the various exponential growth rates compare to each other.  In
other words, one wants to know how common an occurence it is, as a
function of $k$, for an element of $\Mod_g$ to act trivially on the
first $k$ terms of the lower central series of $\pi_1\Sigma_g$.  

\begin{problem}[Logarithmic densities of the Johnson filtration]
Determine the asymptotics of $d_{\log}(\T_g(m))$ both as $g\to \infty$
and as $m\to\infty$.
\end{problem}

Indeed, as far as I know, even the asymptotics of the 
(logarithmic) density of the $k$th term
of the lower central series of $\pi_1\Sigma_g$ in $\pi_1\Sigma_g$ as
$k\to\infty$ has not been determined.

\medskip
\noindent
{\bf Entropy. }
The {\em exponential growth rate} of a 
group $\Gamma$ with respect to a finite generating set $S$ is defined 
as 
$$w(\Gamma,S):=\lim_{r\to\infty}(B_r(\Gamma,S))^{1/r}$$ 
where $B_r(\Gamma,S)$ denotes the cardinality of the $r$-ball in the
Cayley graph of $\Gamma$ with respect to the generating set $S$; the
limit exists since $\beta$ is submultiplicative.  The {\em entropy} of
$\Gamma$ is defined to be 
$$\ent(\Gamma)=\inf\{\log w(\Gamma,S): \mbox{$S$ is finite and generates $\Gamma$}\}$$
 
Among other things, the group-theoretic entropy of 
$\ent(\pi_1M)$ of a closed, Riemannian manifold $M$ gives a lower bound 
for (the product of diameter times) both the volume growth entropy of
$M$ and the topological entropy of the geodesic flow on $M$.  See 
\cite{Harp1} for a survey.  

Eskin-Mozes-Oh \cite{EMO} 
proved that nonsolvable, finitely-generated linear groups
$\Gamma$ have positive entropy.  Since it is classical that 
the action of $\Mod_g$ on $H_1(\Sigma_g,\Z)$ gives a surjection 
$\Mod_g\to \Sp(2g,\Z)$, it follows immediately 
that $\ent(\Mod_g)>0$.  This method of proving positivity of entropy 
fails for the Torelli group $\T_g$ since it is in the kernel of the
standard symplectic representation of $\Mod_g$.  However, one can
consider the action of $\T_g$ on the homology of the universal 
abelian cover of $\Sigma_g$, considered as a (finitely generated) module
over the corresponding covering group, to find a linear representation
of $\T_g$ which is not virtually solvable.  This is basically the {\em
Magnus representation}.  Again by Eskin-Mozes-Oh we
conclude that $\ent(\T_g)>0$.  

\begin{problem}
\label{problem:entropy}
Give explicit upper and lower bounds for $\ent(\Mod_g)$.  Compute the
asymptotics of $\ent(\Mod_g)$ and of $\ent(\T_g)$ as $g\to \infty$.  
Similarly for $\ent(\T_g(k))$ as $k\to\infty$.
\end{problem}

\section{Problems on the geometry of $\M_g$}
\label{section:geomtop}

It is a basic question to understand properties of complex
analytic and geometric structures on $\M_g$, and how these structures
constrain, and are constrained by, the global topology of $\M_g$.  Such
structures arise frequently in applications.  For example one first 
tries to put a geometric structure on $\M_g$, such as that of 
a complex orbifold or of a negatively curved Riemannian manifold.  Once
this is done, general theory for such structures (e.g. Schwarz Lemmas
or fixed point theorems) can then be applied to prove hard theorems.  
Arakelov-Parshin Rigidity and Nielsen Realization are two examples of
this.  In this section we pose a few problems about the topology and
geometry of $\M_g$.

\subsection{Isometries}

Royden's Theorem \cite{Ro} states that when $g\geq 3$, 
every isometry of Teichm\"{u}ller space $\Teich_g$, 
endowed with the Teichm\"{u}ller metric $d_{\Teich}$, is
induced by an element of $\Mod_g$; that is
$$\Isom(\Teich_g,d_{\Teich})=\Mod_g$$

Note that $d_{\Teich}$ comes from a non-Riemannian Finsler metric,
namely a norm on the cotangent space at each point $X\in\Teich_g$.  
This cotangent space can be identified with the space
$Q(X)$ of holomorphic quadratic differentials on $X$.  

I believe Royden's theorem can be generalized from the Teichm\"{u}ller
metric to all metrics.

\begin{conjecture}[Inhomogeneity of all metrics]
\label{conjecture:royden}
Let $\Teich_g$ denote the Teichm\"{u}ller space of closed, 
genus $g\geq 2$ Riemann
surfaces.  Let $h$ be any Riemannian metric (or any 
Finsler metric with some weak regularity 
conditions) on $\Teich_g$ which is invariant under 
the action of the mapping class group $\Mod_g$, and for which this
action has finite covolume.  Then $\Isom(\Teich_g,h)$ is
discrete; even better, it contains $\Mod_g$ as a subgroup of index $C=C(g)$.  
\end{conjecture}

Royden's Theorem is the special case when $h$ is the Teichm\"{u}ller
metric (Royden gets $C=2$ here).  A key philosophical implication 
of Conjecture \ref{conjecture:royden} is that the mechanism behind 
the inhomogeneity of Teichm\"{u}ller space is due not to fine regularity
properties of the unit ball in $Q(X)$ (as
Royden's proof suggests), but to the global topology of moduli space.  This in 
turn is tightly controlled by the structure of $\Mod_g$.  
As one piece of evidence for Conjecture
\ref{conjecture:royden}, I would like to point out that it would follow
if one could extend the main theorem of \cite{FW1} from the closed to the
finite volume case.  

In some sense looking at $\Mod_g$-invariant metrics seems too strong,
especially since $\Mod_g$ has torsion. Perhaps, for example, 
the inhomogeneity of $\Teich_g$ is
simply caused by the constraints of the torsion in $\Mod_g$.
Sufficiently large index subgroups of $\Mod_g$ are torsion free.  Thus
one really wants to strengthen Conjecture \ref{conjecture:royden} by
replacing $\Mod_g$ by any finite index subgroup $H$, and by replacing
the constant $C=C(g)$ by a constant $C=C(g,[\Mod_g:H])$.  After this one
can explore metrics invariant by much smaller subgroups, such as $\T_g$,
and at least hope for discreteness of the corresponding isometry group
(as long as the subgroup is sufficiently large).

If one can prove the part of Conjecture \ref{conjecture:royden} which
gives discreteness of the isometry group of any $\Mod_g$-invariant
metric on $\Teich_g$,  one can approach the stronger statement that
$C_g=1$ or $C_g=2$ as follows.  Take the quotient of $\Teich_g$ by any
group $\Lambda$ of isometries properly containing $\Mod_g$.  
By discreteness of $\Lambda$, the quotient $\Teich_g/\Lambda$ is 
a smooth orbifold which is finitely orbifold-covered by $\M_g$.  

\begin{conjecture}[$\M_g$ is maximal]
For $g\geq 3$ the smooth orbifold $\M_g$ does not finitely 
orbifold-cover any other smooth orbifold.
\end{conjecture}

A much stronger statement, which may be true, would be to prove that if 
$N$ is any finite cover of $\M_g$, then the only orbifolds which $N$ can
orbifold cover are just the covers of $\M_g$.  
Here is a related basic topology question about $\M_g$.

\begin{question}
Let $Y$ be any finite cover of $\M_g$, and let $f:Y\to Y$ be a finite
order homeomorphism.  If $f$ is homotopic to the identity, must $f$
equal the identity?
\end{question}

\subsection{Curvature and $\Q$-rank}

\noindent
{\bf Nonpositive sectional curvature. }There has been a long history of
studying metrics and curvature 
on $\M_g$ \footnote{As 
$\M_g$ is an orbifold, 
technically one studies $\Mod_g$-invariant metrics on 
the Teichm\"{u}ller space $\Teich_g$.}; see, e.g.,  
\cite{BrF,LSY,Mc2}.  A recurring theme is to find aspects of negative
curvature in $\M_g$.  While $\M_g$ admits no metrics of negative
curvature, even in a coarse sense, if $g\geq 2$ (see \cite{BrF}), 
the following remains a basic open problem.

\begin{conjecture}[Nonpositive curvature]
\label{conjecture:npc}
For $g\geq 2$ the orbifold $\M_g$ admits no complete, finite volume
Riemannian metric with nonpositive sectional curvatures uniformly 
bounded away from $-\infty$.  
\end{conjecture}

One might be more ambitious in stating Conjecture \ref{conjecture:npc}, 
by weakening the finite volume condition, by dropping the uniformity of the
curvature bound, or by extending the statement from $\M_g$ to any finite
cover of $\M_g$ (or perhaps even to certain infinite covers).  It would
also be interesting to extend Conjecture \ref{conjecture:npc} beyond the
Riemannian realm to that of ${\rm CAT}(0)$ metrics; see, e.g.,  
\cite{BrF} for a notion of finite volume which extends to this context.

In the end, it seems that we will have to make do with various relative
notions of nonpositive or negative curvature, as in \cite{MM1,MM2}, or 
with various weaker notions of nonpositive curvature, such as holomorphic,
Ricci, or highly singular versions (see, e.g.,\cite{LSY}), or
isoperimetric type versions such a Kobayashi or Kahler hyperbolicity (see
\cite{Mc2}).  Part of the difficulty with trying to fit $\M_g$ into
the ``standard models'' seems to be the topological structure of the 
cusp of $\M_g$.

\bigskip
\noindent
{\bf Scalar curvature and $\Q$-rank. }
Let $S$ be a genus $g$ surface with $n$ punctures, and let $\M(S)$
denote the corresponding moduli space.  We set 
$$d(S) = 3g-3+n$$

The constant $d(S)$ is fundamental
in Teichm\"uller theory: it is the complex dimension of the
Teichm\"uller space $\Teich(S)$; it is also the number of 
curves in any pair-of-pants decomposition of $S$.  While
previous results have concentrated on sectional and holomorphic
curvatures, Shmuel
Weinberger and I have recently proven the following; see \cite{FW2}.

\begin{theorem}[Positive scalar curvature]
\label{conj:psc}
Let $M$ be any finite cover of $\M(S)$.  Then $M$ admits a complete, 
finite-volume Riemannian metric of (uniformly bounded) 
positive scalar curvature if and only if $d(S)\geq 3$. 
\end{theorem}

The analogous statement was proven for locally symmetric 
arithmetic manifolds $\Gamma\backslash G/K$ by Block-Weinberger
\cite{BW}, where $d(S)$ is replaced by the $\Q$-rank of $\Gamma$.  When
$d(S)\geq 3$, the metric on $M$ has the quasi-isometry type of a ray, so
that it is not quasi-isometric to the Teichm\"{u}ller metric on $\M_g$.
It seems likely that this is not an accident.

\begin{conjecture}
\label{conjecture:chang}
Let $S$ be any surface with $d(S)\geq 1$.  Then $M$ does not 
admit a finite volume Riemannian metric of (uniformly bounded) positive
scalar curvature in the quasi-isometry class of the Teichm\"{u}ller metric.
\end{conjecture}

The analogue  of Conjecture \ref{conjecture:chang} in the 
$\Gamma\backslash G/K$ case was proven by S. Chang in \cite{Ch}.  
The same method of proof as in \cite{Ch} should 
reduce Conjecture \ref{conjecture:chang} to the following discussion, 
which seems to be of independent interest, and which came out of
discussions with H. Masur.

What does $\M_g$, endowed with the Teichm\"{u}ller metric $d_{\Teich}$, 
look like from far away?  This can be
formalized by Gromov's notion of 
{\em tangent cone at infinity}:
\begin{equation}
\label{eq:cone}
\cone(\M_g):=\lim_{n\to \infty}(\MS,\frac{1}{n}d_{\Teich})
\end{equation}
where the limit is taken in the sense of Gromov-Hausdorff convergence of
pointed metric spaces; here we have fixed a basepoint in $\M_g$ 
once and for all.  This limit is easily shown to make sense and 
exist in our context.  To state our conjectural answer as to what
$\cone(\M_g)$ looks like, we will need the complex of curves on $\Sigma_g$.  
Recall that the {\em complex of curves} 
${\cal C}_g$ for $g\geq 2$ is the simplicial complex with one
vertex for each nontrivial, nonperipheral isotopy class of simple closed
curves on $\Sigma_g$, and with 
a $k$-simplex for every $(k+1)$-tuple of such isotopy classes
for which there are mutually disjoint representatives.  
Note that $\Mod_g$ acts by simplicial automorphisms on 
${\cal C}_g$.  

\begin{conjecture}[\boldmath$\Q$-rank of moduli space]
\label{conj:cone}
$\cone(\M_g)$ is homeomorphic to the (open) cone on the
quotient ${\cal C}_g/\Mod_g$ \footnote{Note: This 
statement is a slight cheat;
the actual version requires the language of orbi-complexes.}.
\end{conjecture}

One can pose a stronger version of Conjecture \ref{conj:cone} that predicts
the precise bilipschitz type of the natural metric on $\cone(\M_g)$; an 
analogous statement for quotients $\Gamma\backslash G/K$ 
of symmetric spaces by 
lattices was proven by Hattori \cite{Hat}.  H. Masur and I have 
identified the right candidate for a coarse fundamental domain needed 
to prove Conjecture \ref{conj:cone}; its description involves certain
length inequalities analogous to those on roots defining Weyl chambers.
Further, the (conjectured) dimensions of the corresponding tangent cones
are $\Q$-rank$(\Gamma)$ and $d(S)$, respectively.  Thus we propose the 
following additions to the list of analogies between arithmetic lattices
and $\Mod_g$.

\smallskip
\begin{center}
\begin{tabular}{c|c}
arithmetic lattices& $\Mod_g$\\
\hline  
$\Q$-rank($\Gamma$) & $d(S)$\\
root lattice & \{simple closed curves\} \\
simple roots & top. types of simple closed curves\\
$\cone(\Gamma\backslash G/K)$& $\cone(\MS)$\\
\end{tabular}
\end{center}

\smallskip
As alluded to above, Conjecture \ref{conj:cone} should imply, 
together with the methods in \cite{Ch}, the second statement of Conjecture
\ref{conj:psc}.

\subsection{The Kahler group problem}

Complete Kahler metrics on $\M_g$ 
with finite volume and bounded curvatures have been found by Cheng-Yau,
McMullen and others (see, e.g., \cite{LSY} for a survey).  
The following conjecture, however, is still not
known.  I believe this is a folklore conjecture.

\begin{conjecture}[$\Mod_g$ is Kahler]
\label{conjecture:kahler}
For $g\geq 3$, the group 
$\Mod_g$ is a {\em Kahler group}, i.e. it is isomorphic to 
the fundamental group of a compact Kahler manifold. 
\end{conjecture}

It was shown in \cite{Ve} that $\Mod_2$ is not a Kahler group.  This
is proven by reducing (via finite extensions) to the pure 
braid group case; these groups are not Kahler since they are iterated
extensions of free groups by free groups.  

A natural place to begin proving Conjecture \ref{conjecture:kahler} is 
the same strategy that Toledo uses in \cite{To} for nonuniform lattices
in ${\rm SU}(n,1), n\geq 3$.  The main point is the following.  One
starts with a smooth open variety $V$ and wants to prove that 
$\pi_1V$ is a Kahler group.  The first step is to find 
a compactification $\bar{V}$ of $V$ which is projective, and for which 
$\bar{V}-V$ has codimension {\em at least $3$}.  This assumption
guarantees that the intersection of the generic $2$-plane $P$ in projective
space with $\bar{V}$ misses $V$.  The (weak) Lefschetz Theorem then
implies that the inclusion $i:\bar{V}\cap P\hookrightarrow V$ induces an
isomorphism on fundamental groups, thus giving that $\pi_1V$ is a
Kahler group.  

One wants to apply this idea to the Deligne-Mumford compactification 
$\overline{\M_g}$ of moduli space $\M_g$.  This almost works, except
that there is a (complex) codimension {\em one} singular stratum of
$\overline{\M_g}$, so that the above does not apply.  Other 
compactifications of $\M_g$ are also problematic in this regard.

What about the Torelli group $\T_g$?  This group, at least for $g\geq
6$, is not known to violate any of the known constraints on Kahler
groups.  Most notably, Hain \cite{Ha3} proved the deep result that 
for $g\geq 6$ the group $\T_g$ has a quadratically
presented Malcev Lie algebra; this is one of the more subtle properties 
posessed by Kahler groups.  Note Akita's theorem (Theorem
\ref{theorem:akita} above) that the classifying space of $\T_g, g\geq 7$
does not have the homotopy type of a finite complex shows that these
$\T_g$ are not fundamental groups of closed {\em aspherical} manifolds.
There are of course Kahler groups (e.g. finite Kahler groups) 
with this property.  In contrast to 
Conjecture \ref{conjecture:kahler}, I have recently 
proven \cite{Fa3} the following.

\begin{theorem}
$\T_g$ is not a Kahler group for any $g\geq 2$.  
\end{theorem}

Denoting the symplectic representation by $\pi:\Mod_g\to
\Sp(2g,\Z)$, the next question is to ask which of the groups
$\pi^{-1}(\Sp(2k,\Z))$ for $1\leq k\leq 2g$ interpolating between the two
extremes $\T_g$ and $\Mod_g$ are Kahler groups.  My only
guess is that when $k=1$ the group is not Kahler.

\subsection{The period mapping}
\label{subsection:period:mapping}

To every Riemann surface $X\in \Teich_g$ one attaches its {\em
Jacobian} $\Jac(X)$ , which is the quotient of the dual
$(\Omega^1(X))^\ast\approx\C^g$ 
of the space of holomorphic $1$-forms on $X$ by the
lattice $H_1(X,\Z)$, where $\gamma\in H_1(X,\Z)$ is the linear functional
on $\Omega^1(X)$ given by $\omega\mapsto\int_\gamma\omega$.  Now
$\Jac(X)$ is a complex torus, and Riemann's period relations show that
$\Jac(X)$ is also an algebraic variety, i.e. $\Jac(X)$ is an {\em abelian
variety}.  The algebraic intersection number on $H_1(X,\Z)$ induces a
symplectic form on $\Jac(X)$, which can be thought of as the imaginary
part of a positive definite Hermitian form on $\C^g$.  This extra bit of
structure is called a {\em principal polarization} on the abelian
variety $\Jac(X)$.  

The space $\A_g$ of all $g$-dimensional (over $\C$) principally
polarized abelian varieties is parameterized by the locally symmetric
orbifold $\Sp(2g,\Z)\backslash \Sp(2g,\R)/{\rm U}(g)$.  The {\em
Schottky problem}, one of the central classical problems of algebraic
geometry, asks for the image of the {\em period mapping}
$$\Psi:\M_g\to\A_g$$ which sends a surface $X$ to its Jacobian
$\Jac(X)$.  In other words, the Schottky problem asks: which principally
polarized abelian varieties occur as Jacobians of some Riemann surface?
Torelli's Theorem states that $\Psi$ is injective; the 
image $\Psi(\M_g)$ is called the {\em period locus}.  
The literature on this problem is vast (see, e.g., \cite{D} for a
survey), and goes well beyond the scope of the present paper.

Inspired by the beautiful paper \cite{BS} of Buser-Sarnak, I
would like to pose some questions about the period locus
from a different (and probably nonstandard) point of view.
Instead of looking for precise algebraic equations describing
$\Psi(\M_g)$, what if we instead try to figure out how to tell
whether or not it contains a given torus, or if we
try to describe what the period locus {\em roughly} looks like? 
Let's make these questions more precise.

The data determining a principally polarized abelian variety is not
combinatorial, but is a matrix of real numbers.  However, one can still
ask for algorithms involving such data by using the {\em complexity theory
over $\R$} developed by Blum-Shub-Smale; see, e.g., \cite{BCSS}.  Unlike
classical complexity theory, here one assumes that real numbers can be
given and computed precisely, and develops algorithms, measures of
complexity, and the whole theory under this assumption.  In the language
of this theory we can then pose the following.

\begin{problem}[Algorithmic Schottky problem]
Give an algorithm, in the
sense of complexity theory over $\R$, 
which takes as input a $2g\times 2g$ symplectic
matrix representing a principally polarized abelian variety, 
and as output tells whether or
not that torus lies in the period locus. 
\end{problem}

One might also fix some $\epsilon=\epsilon(g)$, and then ask
whether or not a given principally polarized abelian variety lies within
$\epsilon$ (in the locally symmetric metric on ${\cal A}_g$) of the
period locus.  

It should be noted that 
S. Grushevsky \cite{Gr} has made the KP-equations solution to the
Schottky problem effective in an algebraic sense.  This seems to be
different than what we have just discussed, though.

We now address the question of what the Schottky locus looks like from
far away.  To make this precise, let $\cone(\A_g)$ denote the tangent 
cone at infinity (defined in
(\ref{eq:cone}) above) of the locally symmetric Riemannian orbifold 
$\A_g$.  Hattori \cite{Hat} proved that $\cone(\A_g)$ is
homeomorphic to the open cone on the quotient of the Tits boundary of 
the symmetric space $\Sp(2g,\R)/{\rm U}(g)$; indeed it is isometric to 
a Weyl chamber in the symmetric space,
which is just a Euclidean sector of dimension $g$.  

\begin{problem}[Coarse Schottky problem]
Describe, as a subset of a $g$-dimensional Euclidean sector, the subset 
of $\cone(\A_g)$ determined by the 
Schottky locus in $\A_g$.
\end{problem}

Points in $\cone(\A_g)$ are recording how the relative sizes of basis
vectors of the tori are changing; it is precisely the ``skewing
parameters'' that are being thrown away.  It doesn't seem unreasonable
to think that much of the complexity in describing the Schottky locus is
coming precisely from these skewing parameters, so that this
coarsification of the Schottky problem, unlike the classical version, 
may have a reasonably explicit solution.

There is a well-known feeling that the
Schottky locus is quite distorted in $\A_g$.  Hain and Toledo (perhaps
among others) have posed the problem of determining the second
fundamental form of the Schottky locus, although they indicate that this
would be a rather difficult computation. We can coarsify this question 
by extending the definition of distortion of subgroups given in 
Subsection \ref{subsection:distortion} above to the context of subspaces
of metric spaces.  Here the distortion of a subset $S$ in a metric space
$Y$ is defined by comparing the restriction of the 
metric $d_Y$ to $S$ versus the 
induced path metric on $S$.  

\begin{problem}[Distortion of the Schottky locus]
Compute the distortion of the Schottky locus in $\A_g$.
\end{problem}

A naive guess might be that it is exponential.
 
\section{The Torelli group}
\label{section:torelli}

Problems about the Torelli group $\T_g$ have a special flavor of their
own.  As one passes from $\Mod_g$ to $\T_g$, significant and beautiful
new phenomena occur.  One reason for the richness of this theory is that
the standard exact sequence $$1\to \T_g\to \Mod_g\to \Sp(2g,\Z)\to 1$$
gives an action $\psi:\Sp(2g,\Z)\to\Out(\T_g)$, so that any natural
invariant attached to $\T_g$ comes equipped with an $\Sp(2g,\Z)$-action.
The most notable examples of this are the cohomology algebra
$H^\ast(\T_g,\Q)$ and the Malcev Lie algebra ${\cal L}(\T_g)\otimes \Q$,
both of which become $\Sp(2g,\Q)$-modules, allowing for the application
of symplectic representation theory.  See, e.g., \cite{Jo1,Ha3,Mo4} for
more detailed explanations and examples.  

In this section I present a few of my favorite problems.  I refer the
reader to the work of Johnson, Hain and Morita for other problems about 
$\T_g$; see, e.g., \cite{Mo1,Mo3}.

\subsection{Finite generation problems}

For some time it was not known
if the group $\K_g$ generated by Dehn twists about bounding curves 
was equal to, or perhaps a finite index subgroup of $\T_g$, 
until Johnson found the {\em Johnson homomorphism} $\tau$ and
proved exactness of:
\begin{equation}
\label{eq:johnson:exact}
1\to\kg\to\T_g\stackrel{\tau}{\to}\wedge^3H/H \to 1
\end{equation}
where $H=H_1(\Sigma_g;\Z)$ and where the inclusion $H\hookrightarrow
\wedge^3H$ is given by the map $x\mapsto x\wedge \hat{i}$, where
$\hat{i}$ is the intersection from 
$\hat{i}\in \wedge^2H$.  Recall that a 
{\em bounding pair map} is a composition 
$T_a\circ T_b^{-1}$ of Dehn twists about {\em bounding pairs}, i.e. 
pairs of disjoint, nonseparating, 
homologous simple closed curves $\{a,b\}$.  Such a
homeomorphism clearly lies in $\T_g$.  By direct
computation Johnson shows that the $\tau$-image of such a map is
nonzero, while $\tau(\K_g)=0$; proving that $\ker(\tau)=\K_g$ is much
harder to prove.

Johnson proved in \cite{Jo2} that
$\T_g$ is finitely generated for all $g\geq 3$.  McCullough-Miller 
\cite{McM} proved that $\K_2$ is not
finitely generated; Mess \cite{Me} then proved that $\K_2$ is in fact a
free group of infinite rank, generated by the set of symplectic
splittings of $H_1(\Sigma_2,\Z)$.  The problem of finite 
generation of $\K_g$ for all $g\geq 3$ was 
recently solved by Daniel Biss and me in \cite{BF}.

\begin{theorem}
\label{theorem:biss}
The group $\K_g$ is not finitely generated for any $g\geq 2$.  
\end{theorem}

The following basic problem, however, remains open (see, e.g., \cite{Mo3}, 
Problem 2.2(ii)).

\begin{question}[Morita]
\label{q:h1}
Is $H_1(\K_g,\Z)$ finitely generated for $g\geq 3$?
\end{question}

Note that Birman-Craggs-Johnson (see, e.g., \cite{BC,Jo1}) and 
Morita \cite{Mo2} have found large abelian quotients of 
$\K_g$. 

The proof in \cite{BF} of Theorem \ref{theorem:biss}, suggests an approach to 
answering to Question \ref{q:h1}.  Let me briefly describe the 
idea.  Following the the outline 
in \cite{McM}, we first find an action of $\K_g$ on the first homology
of a certain abelian cover $Y$ of $\Sigma_g$; this action respects 
the structure of $H_1(Y,\Z)$ as a module over the 
Galois group of the cover.  The crucial piece is that we are 
able to reduce this to a representation
$$\rho:\K_g\to \SL_2(\Z[t,t^{-1}])$$
on the special linear group over the ring of integral laurent
polynomials in one variable.  This group acts on an associated
Bruhat-Tits-Serre tree, and one can then analyze this action using 
combinatorial group theory.  

One might now try to answer Question \ref{q:h1} in the negative 
by systematically computing
more elements in the image of $\rho$, and then analyzing more closely 
the action on the tree for $\SL_2$.  One potentially useful ingredient 
is a theorem of Grunewald-Mennike-Vaserstein \cite{GMV} which gives free
quotients of arbitrarily high rank for the group $\SL_2(\Z[t])$ and the 
group $\SL_2(K[s,t])$, where $K$ is
an arbitrary finite field.

Since we know for $g\geq 3$ that $\T_g$ is finitely generated and
$\K_g$ is not, it is natural to ask about the subgroups interpolating
between these two. To be precise, consider the exact sequence
(\ref{eq:johnson:exact}).  Corresponding to each 
subgroup $L<\wedge^3H/H$ is its pullback $\pi^{-1}(L)$.  The
lattice of such subgroups $L$ can be thought of as a kind of
interpolation between $\T_g$ and $\K_g$.  

\begin{problem}[Interpolations]
Let $g\geq 3$. For each subgroup $L<\wedge^3H/H$, 
determine whether or not $\pi^{-1}(L)$ is finitely generated.
\end{problem}

As for subgroups deeper down than $\K_g=\T_g(2)$ 
in the Johnson filtration $\{\T_g(k)\}$, we would like to record the following.

\begin{theorem}[Johnson filtration not finitely generated]
\label{theorem:jf:notfg}
For each $g\geq 3$ and each $k\geq 2$, the group $\T_g(k)$ is not
finitely generated.
\end{theorem}

\proof
We proceed by induction on $k$.  For $k=2$ this is just the theorem of
\cite{BF} that $\K_g=\T_g(2)$ is not finitely generated for any $g\geq
3$.  Now assume the theorem is true for $\T_g(k)$.  
The {\em $k$th Johnson homomorphism} is a homomorphism 
$$\tau_g(k):\T_g(k)\to {\mathfrak h}_g(k)$$
where ${\mathfrak h}_g(k)$ is a certain finitely-generated abelian group,
coming from the $k^{th}$ graded piece of a certain graded Lie algebra;
for the precise definitions and constructions, see, e.g., \S 5 of \cite{Mo3}.  
All we will need is Morita's result (again, see \S 5 of \cite{Mo3}) that 
$\ker(\tau_g(k))=\T_g(k+1)$.  We thus have an exact sequence
\begin{equation}
\label{eq:bf2}
1\to\T_g(k+1)\to\T_g(k)\to\tau_g(k)(\T_g(k))\to 1
\end{equation}
Now the image $\tau_g(k)(\T_g(k))$ is a subgroup of the finitely generated
abelian group ${\mathfrak h}_g(k)$, and so is finitely generated.  
If $\T_g(k+1)$ were finitely-generated, then by (\ref{eq:bf2}) we would
have that $\T_g(k)$ is finitely generated, contradicting the induction
hypothesis.  Hence $\T_g(k+1)$ cannot be finitely generated, and we are
done by induction.
\endproof

The Johnson filtration $\{\T_g(k)\}$ and 
the lower central series $\{(\T_g)_k\}$ do not coincide; indeed Hain
proved in \cite{Ha3} that there are terms of the former not contained in
any term of the latter.  Thus the following remains open.

\begin{conjecture}
\label{conjecture:fglcc}
For each $k\geq 1$, the group $(\T_g)_k$ is not finitely generated.
\end{conjecture}

Another test of our understanding of the Johnson filtration is the
following. 

\begin{problem}
Find $H_1(\T_g(k),\Z)$ for all $k\geq 2$.
\end{problem}

\medskip
\noindent
{\bf Generating sets for $\T_g$. }
One difficulty in working with $\T_g$ is the complexity of its
generating sets: any such set must have at least $\frac{1}{3}[4g^3-g]$ elements
since $\T_g$ has abelian quotients of this rank 
(see \cite{Jo5}, Corollary after Theorem 5).  Compare this with
$\Mod_g$, which can always be generated by $2g+1$ Dehn twists
(Humphries), or even by $2$ elements (Wajnryb \cite{Wa})!  
How does one keep track, for example, of the (at least) $1330$
generators for $\T_{10}$?  How does one even give a usable naming 
scheme for working
with these?  Even worse, in Johnson's proof of finite generation of
$\T_g$ (see \cite{Jo2}), the given generating set has $O(2^g)$
elements.  The following therefore seems fundamental; at the very least
it seems that solving it will require us to understand the
combinatorial topology underlying $\T_g$ in a deeper way than we now
understand it.

\begin{problem}[Cubic genset problem]
\label{p:gen}
Find a generating set for $\T_g$ with $O(g^d)$ many elements for some $d\geq
3$. Optimally one would like $d=3$.
\end{problem}

In fact in \S 5 of \cite{Jo2}, Johnson explicitly poses a much harder
problem: for $g\geq 4$ can $\T_g$ be generated by $\frac{1}{3}[4g^3-g]$
elements?  As noted above, this would be a sharp result.  Johnson
actually obtains this sharp result in genus three, by finding
(\cite{Jo2}, Theorem 3) an explicit set of $35$ generators for $\T_g$.
His method of converting his $O(2^g)$ generators to $O(g^3)$ becomes far
too unwieldy when $g>3$.

One approach to Problem \ref{p:gen} is to 
follow the original plan of \cite{Jo2}, but using a simpler generating
set for $\Mod_g$.  This was indeed the motivation for Brendle and me when
we found in \cite{BFa1} a generating set for $\Mod_g$ consisting of $6$ 
involutions, i.e. $6$ elements of order $2$.  This bound was later
improved by Kassabov \cite{Ka} to $4$ elements of
order $2$, at least when $g\geq 7$.  Clearly $\Mod_g$ is never
generated by $2$ elements of order two, for then it would be a
quotient of the infinite dihedral group, and so would be virtually
abelian.  Since the current known bounds are so close to being sharp, 
it is natural to ask for the sharpest bounds.

\begin{problem}[Sharp bounds for involution generating sets]
For each $g\geq 2$, prove sharp bounds for the minimal number of
involutions required to generate $\Mod_g$.  In particular, for $g\geq
7$ determine whether or not $\Mod_g$ is generated by $3$ involutions.
\end{problem}

\subsection{Higher finiteness properties and cohomology}

While there has been spectacular progress in understanding
$H^\ast(\Mod_g,\Z)$ (see \cite{MW}), very little is known about 
$H^\ast(\T_g,\Z)$, and even less about $H^\ast(\K_g,\Z)$.  Further, we
do not have answers to the basic finiteness questions one asks about 
groups.  

Recall that the {\em cohomological dimension} of
a group $\Gamma$, denoted $\cd(\Gamma)$, is defined to be
$$\cd(\Gamma):=\sup\{i:H^i(\Gamma,V)\neq 0 \mbox{\ for some
$\Gamma$-module $V$}\}$$
If a group $\Gamma$ has a torsion-free subgroup $H$ of finite index, then
the {\em virtual} cohomological dimension of $\Gamma$ is defined to be
$\cd(H)$; Serre proved that this number does not depend on the choice of
$H$.  It is a theorem of Harer, with earlier estimates and a later different
proof due to Ivanov, that $\Mod_g$ has virtual cohomological dimension
$4g-5$; see \cite{Iv1} for a summary.

\begin{problem}[Cohomological Dimension]
\label{problem:cd}
Compute the cohomological dimension of $\T_g$ and of $\K_g$.  More
generally, compute the cohomological dimension of $\T_g(k)$ for all 
$k\geq 1$.
\end{problem}

Note that the cohomological dimension $\cd(\T_g)$ is bounded above by
the (virtual) cohomological dimension of $\Mod_g$, which is $4g-5$.  
The following is a start on some lower bounds.  

\begin{theorem}[Lower bounds on $\cd$]
For all $g\geq 2$, the following inequalities hold:
\begin{displaymath}
\begin{array}{l}
{\cd(\T_g)\geq \left\{\begin{array}{ll}
(5g-8)/2& \textrm{if $g$ is even}\\
&\\
(5g-9)/2&\textrm{if $g$ is odd}
\end{array}\right.}\\
\\
\cd(\K_g)\geq 2g-3
\\
\\
\cd(\T_g(k))\geq g-1 \textrm{\ \ for $k\geq 3$}
\end{array}
\end{displaymath}

\end{theorem}

\proof
Since for any group $\Gamma$ with $\cd(\Gamma)<\infty$ we have
$\cd(\Gamma)\geq \cd(H)$ for any subgroup $H<\Gamma$, an easy way 
to obtain lower bounds for $\cd(\Gamma)$ is to find large
free abelian subgroups of $\Gamma$.  To construct such subgroups for
$\T_g$ and for $\K_g$, take a maximal collection of mutually disjoint 
separting curves on $\Sigma_g$; by an Euler characteristic argument it
is easy to see that there are $2g-3$ of these, and it is not hard to
find them.  The group generated by Dehn twists on these curves is
isomorphic to $\Z^{2g-3}$, and is contained in $\K_g<\T_g$.  

For $\T_g$ we obtain the better bounds by giving a slight variation of 
Ivanov's discussion of {\em Mess subgroups}, given in \S 6.3 of
\cite{Iv1}, adapted so that the constructed subgroups lie in $\T_g$.    
Let $\Mod^1_g$ denote the
group of homotopy classes of orientation-preserving homeomorphisms of
the surface $\Sigma^1_g$ of genus $g$ with one boundary component, 
fixing $\partial \Sigma^1_g$ pointwise, up to isotopies which fix
$\partial \Sigma^1_g$ pointwise.  
We then have a well-known exact sequence (see,
e.g. \cite{Iv1}, \S 6.3)
\begin{equation}
\label{eq:bndry}
1\to\pi_1T^1\Sigma_g\to\Mod^1_g\stackrel{\displaystyle
\pi}{\displaystyle \to}\Mod_g\to 1
\end{equation}
where $T^1\Sigma_g$ is the unit tangent bundle of $\Sigma_g$.  
Now suppose $g\geq 2$.  Let $C_2$ and $C_3$ be maximal abelian subgroups
of $\T_2$ and $\T_3$, respectively; these have ranks $1$ and $3$,
respectively.  We now define $C_g$ inductively, beginning with $C_2$ if
$g$ is even, and with $C_3$ if $g$ is odd.  
Let $C^1_g$ be the pullback $\pi^{-1}(C_g)$ of $C_g$ 
under the map $\pi$ in (\ref{eq:bndry}).  Note that, since the copy of
$\pi_1T^1\Sigma_g$ in $\Mod^1_g$ is generated by ``point pushing'' and
the twist about $\partial \Sigma_g$, it actually lies in the
corresponding Torelli group $\T^1_g$.  The inclusion
$\Sigma^1_g\hookrightarrow \Sigma_{g+2}$ induces an injective
homomorphism $i:\T^1_g\hookrightarrow \T_{g+2}$ 
via ``extend by the identity''. The
complement of $\Sigma^1_g$ in $\Sigma_{g+2}$ clearly contains a pair of
disjoint separating curves.  Now define $C_{g+2}$ to be the group by the
Dehn twists about these curves together with $i(C^1_g)$.  Thus
$C_{g+2}\approx C^1_g\times \Z^2$.  The same exact argument as in the
proof of Theorem 6.3A in \cite{Iv1} gives the claimed answers for
$\cd(C_g)$.  

Finally, for the groups $\T_g(k)$ with $k,g\geq 3$ we make the following
construction.  $\Sigma_g$ admits a homeomorphism $f$ of order $g-1$, given
by rotation in the picture of a genus one subsurface $V$ with $g-1$ boundary
components, with a torus-with-boundary attached to each component of
$\partial V$.  It is then easy to see that there is a collection of $g-1$
mutually disjoint, $f$-invariant collection of simple closed curves
which decomposes $\Sigma_g$ into a union of $g-1$ subsurfaces
$S_1,\ldots ,S_{g-1}$, each having 
genus one and two boundary components, with mutually disjoint
interiors.  

Each $S_i$ contains a pair of separating curves
$\alpha_i,\beta_i$ with $i(\alpha_i,\beta_i)\geq 2$.  
The group generated by the Dehn twists about
$\alpha_i$ and $\beta_i$ thus generates a free group $L_i$ of rank $2$ (see,
e.g. \cite{FMa}).  Nonabelian free groups have 
elements arbitrarily far down in their lower central series.  As proven 
in Lemma 4.3 of \cite{FLM}, any element in the $k$th level of the 
lower central series for any $L_i$ gives an element $\gamma_i$ lying in 
$\T_g(k)$.  Since $i(\gamma_i,\gamma_j)=0$ for each $i,j$, it follows
that the group $A$ generated by Dehn twists about each $\gamma_i$ is 
isomorphic to $\Z^{g-1}$.  As $A$ can be chosen to lie in 
in any $\T_g(k)$ with $k\geq 3$, we are done.
\endproof

Since $\cap_{k=1}^\infty\T_g(k)=0$, we know
that there exists $K>1$ with the property that the cohomological
dimension of $\T_g(k)$ is constant for all $k\geq K$.  It would be
interesting to determine the smallest such $K$. 
A number of people have different guesses about what the higher
finiteness properties of $\T_g$ should be.  

\begin{problem}[Torelli finiteness]
\label{problem:finiteness1}
Determine the maximal number $f(g)$ for which there is a $K(\T_g,1)$
space with finitely many cells in dimensions $\leq f(g)$.
\end{problem}

Here is what is currently known about Problem \ref{problem:finiteness1}:
\begin{enumerate}
\item $f(2)=0$ since $\T_2$ is not finitely generated (McCullough-Miller
\cite{McM}). 

\item $f(3)\leq 3$ (Johnson-Millson, unpublished, referred to in
\cite{Me}). 

\item For $g\geq 3$, combining Johnson's finite generation result
\cite{Jo2} and a theorem of Akita (Theorem \ref{theorem:akita} below) gives 
$1\leq f(g)\leq 6g-5$.

\end{enumerate}

One natural guess which fits with the (albeit small amount of) 
known data is that $f(g)=g-2$.  
As a special case of Problem \ref{problem:finiteness1}, we emphasize the
following, which is a folklore conjecture.

\begin{conjecture}
$\T_g$ is finitely presented for $g\geq 4$.
\end{conjecture}

One thing we do know is that, in contrast to $\Mod_g$, neither 
$\T_g$ nor $\K_g$ has a classifying space which is homotopy equivalent
to a finite complex; indeed Akita proved the following stronger result.

\begin{theorem}[Akita \cite{Ak}]
\label{theorem:akita}
For each $g\geq 7$ the vector space $H_\ast(\T_g,\Q)$ is infinite
dimensional.  For each $g\geq 2$ the vector space $H_\ast(\K_g,\Q)$ is 
infinite dimensional.
\end{theorem}

Unfortunately the proof of Theorem \ref{theorem:akita} does not
illuminate the reasons why the theorem is true, especially since the
proof is far from constructive.  In order to demonstrate this, and
since the proof idea is simple and pretty, we sketch the proof.

\bigskip
\noindent
{\bf Proof sketch of Theorem \ref{theorem:akita} for \boldmath$\T_g$. }We give
the main ideas of the proof, which is based on a similar argument made
for $\Out(F_n)$ by Smillie-Vogtmann; see \cite{Ak} for details and 
references.
    
If $\dim_{\Q}(H_\ast(\T_g,\Q))<\infty$, then the multiplicativity of 
the Euler characteristic for group extensions, applied to 
$$1\to \T_g\to \Mod_g\to \Sp(2g,\Z)\to 1$$
gives that
 
\begin{equation}
\label{eq:akita1}
\chi(\T_g)=\chi(\Mod_g)/\chi(\Sp(2g,\Z))
\end{equation}

Each of the groups on the right hand side of (\ref{eq:akita1}) has been
computed; the numerator by Harer-Zagier and the denominator by Harder.  
Each of these values is given as a value of the Riemann zeta function
$\zeta$.  Plugging in these values into (\ref{eq:akita1}) gives 
\begin{equation}
\label{eq:akita2}
\chi(\T_g)=\frac{\displaystyle 1}{\displaystyle
2-2g}\prod_{k=1}^{g-1}\frac{\displaystyle 1}{\displaystyle
\zeta(1-2k)}
\end{equation}

It is a classical result of Hurwitz that each finite order element in
$\Mod_g$ acts nontrivially on $H_1(\Sigma_g,\Z)$; hence $\T_g$ is
torsion-free. Thus $\chi(\T_g)$ is an integer.  The rest of the proof of
the theorem consists of using some basic properties of $\zeta$ to prove
that the right hand side of (\ref{eq:akita2}) is not an integer.  
\endproof

The hypothesis $g\geq 7$ in Akita's proof is used only in showing that 
the right hand side of (\ref{eq:akita2}) is not an integer.  This might 
still hold for $g<7$.

\begin{problem}
Extend Akita's result to $2<g<7$.
\end{problem}

Since Akita's proof produces no explicit homology classes, the following seems 
fundamental.

\begin{problem}[Explicit cycles]
Explicitly construct infinitely many linearly
independent cycles in $H_\ast(\T_g,\Q)$ and $H_\ast(\K_g,\Q)$.  
\end{problem}

So, we are still at the stage of trying to find 
explicit nonzero cycles.  In a series of papers (see \cite{Jo1} for a
summary), Johnson proved the quite nontrivial result:
\begin{equation}
\label{eq:h1:torelli}
H^1(\T_g,\Z)\approx \frac{\displaystyle\wedge^3H}{\displaystyle H}\oplus B_2
\end{equation}
where $B_2$ consists of $2$-torsion.  While the $\wedge^3H/H$ piece comes
from purely algebraic considerations, the $B_2$ piece is ``deeper'' in
the sense that it is purely topological, and 
comes from the Rochlin invariant (see \cite{BC} and \cite{Jo1});
indeed the former appears in $H_1$ of the ``Torelli group'' in 
the analogous theory for $\Out(F_n)$, while
the latter does not.

\medskip
\noindent
{\bf Remark on two of Johnson's papers. }
While Johnson's computation of $H_1(\T_g,\Z)$ and his
theorem that $\ker(\tau)=\K_g$ are fundamental results in this area, 
I believe that the details of the proofs of these results are not
well-understood.  These results are proved in \cite{Jo4} and
\cite{Jo3}, respectively; the paper \cite{Jo4} is a particularly dense
and difficult read.  
While Johnson's work is always careful and detailed, and so the
results should be accepted as true, I think it would be
worthwhile to understand \cite{Jo3} and \cite{Jo4}, to exposit them in a
less dense fashion, and perhaps even to give new proofs of their main 
results.  For \cite{Jo3} this is to some extent accomplished
in the thesis \cite{vdB} of van den Berg, where she takes a different approach 
to computing $H_1(\T_g,\Z)$.  

Since dimension one is the only dimension $i\geq 1$ for which we
actually know the $i^{th}$ cohomology of $\T_g$, and since very general
computations seem out of reach at this point, the following seems like a
natural next step in understanding the cohomology of $\T_g$.

\begin{problem}
\label{p:g}
Determine the subalgebras of $H^\ast(\T_g,K)$, for $K=\Q$ and
$K=\F_2$, generated by $H^1(\T_g,K)$.
\end{problem}

Note that $H^\ast(\T_g,K)$ is a module 
over $\Sp(2g,K)$.  When $K=\Q$ this problem has been solved in degree
$2$ by Hain \cite{Ha3} and degree $3$ (up to one unknown piece) 
by Sakasai \cite{Sa}. Symplectic representation theory
(over $\R$) is used as a tool in these papers to greatly simplify 
computations.    When $K=\F_2$, the seemingly basic facts one needs about
representations are either false or they
are beyond the current methods of modular representation theory.  Thus 
computations become more complicated.  Some progress in this case 
is given in \cite{BFa2}, where direct geometric computations,
evaluating cohomology classes on abelian cycles, shows that 
each of the images of 
$$\sigma^\ast:H^2(B_3,\F_2)\to H^2(\T_g,\F_2)$$
$$(\sigma|_{\K_g})^\ast:H^2(B_2,\F_2)\to H^2(\K_g,\F_2)$$
has dimension at least $O(g^4)$.

\subsection{Automorphisms and commensurations of the Johnson
filtration}
The following theorem indicates that all of the algebraic structure of the 
mapping class group $\Mod_g$ is already determined by the infinite index
subgroup $\T_g$, and indeed the infinite index subgroup $\K_g$ of
$\T_g$.  Recall that the {\em extended mapping class group}, denoted
$\Mod_g^\pm$, is defined as the group of homotopy classes of
{\em all} homeomorphisms of $\Sigma_g$, including the
orientation-reversing ones; it contains $\Mod_g$ as a subgroup of index $2$.

\begin{theorem}
\label{thm:fi}
Let $g\geq 4$. Then $\Aut(\T_g)\approx \Mod_g^\pm$ and
$\Aut(\K_g)\approx \Mod_g$.  In fact $\Comm(\T_g)\approx \Mod_g^\pm$ and
$\Comm(\K_g)\approx \Mod_g$.
\end{theorem}

The case of $\T_g, g\geq 5$ was proved by Farb-Ivanov \cite{FI}.
Brendle-Margalit \cite{BM} built on \cite{FI} to prove the harder results on 
$\K_g$.  The cases of $\Aut$, where one can use explicit relations, 
were extended to $g\geq3$ by McCarthy-Vautaw \cite{MV}.  Note too that 
$\Aut(\Mod_g)=\Mod_g^\pm$, as shown by Ivanov (see \S 8 of \cite{Iv1}).

\begin{question}
For which $k\geq 1$ is it true that $\Aut(\T_g(k))=\Mod_g^\pm$? that 
$\Comm(\T_g(k))=\Mod_g^\pm$?
\end{question}

Theorem \ref{thm:fi} answers the question for $k=1,2$.  It would be
remarkable if all of $\Mod_g$ could be reconstructed from
subgroups deeper down in its lower central series.

\subsection{Graded Lie algebras associated to $\T_g$}
\label{section:liealg}

Fix a prime $p\geq 2$.  For a group $\Gamma$ let 
$P_i(\Gamma)$ be defined inductively via $P_1(\Gamma)=\Gamma$ and 
$$P_{i+1}(\Gamma):=[\Gamma,P_i(\Gamma)]\Gamma^p \mbox{\ \ for $i\geq
1$}$$ 
The sequence $\Gamma\supseteq P_2(\Gamma)\supseteq\cdots $
is called the {\em lower exponent $p$ central series}.   The quotient
$\Gamma/P_2(\Gamma)$ has the universal property that any homomorphism
from $\Gamma$ onto an abelian $p$-group factors through
$\Gamma/P_2(\Gamma)$; the group $\Gamma/P_{i+1}(\Gamma)$ has the analgous 
universal property for homomorphisms from $\Gamma$ onto class $i$ nilpotent
$p$-groups.   We can form the direct sum of vector spaces over the
field $\F_p$:
$${\cal L}_p(\Gamma):=\bigoplus_{i=1}^\infty 
\frac{P_i(\Gamma)}{P_{i+1}(\Gamma)}$$
The group commutator on $\Gamma$ induces a bracket on ${\cal L}_p(\Gamma)$
under which it becomes a graded Lie algebra over $\F_p$. See,
e.g. \cite{Se} for the basic theory of Lie algebras over $\F_p$.  

When $p=0$, that is when $P_{i+1}(\Gamma)=[\Gamma,P_i(\Gamma)]$, we
obtain a graded Lie algebra ${\cal L}(\Gamma):={\cal L}_0(\Gamma)$ 
over $\Z$.  The Lie algebra 
${\cal L}(\Gamma)\otimes \R$ is isomorphic to the associated graded
Lie algebra of the Malcev Lie algebra associated to $\Gamma$.    
In \cite{Ha3} Hain found a presentation for the infinite-dimensional 
Lie algebra ${\cal L}(\T_g)\otimes \R$: it is
(at least for $g\geq 6$) 
the quotient of the free Lie algebra on $H_1(\T_g,\R)=\wedge^3H/H, 
H:=H_1(\Sigma_g,\R)$, modulo a finite set of {\em quadratic relations}, 
i.e. modulo an ideal generated by certain elements lying 
in $[P_2(\T_g)/P_3(\T_g)]\otimes \R$.  Each of these relations can already be
seen in the Malcev Lie algebra of the pure braid group.

The main ingredient in Hain's determination of ${\cal L}(\T_g)\otimes
\R$ is to apply Deligne's {\em mixed Hodge theory}.  This theory is a
refinement and extension of the classical Hodge decomposition.  For each
complex algebraic variety $V$ it produces, in a functorial way, various
filtrations with special properties on $H^\ast(V,\Q)$ and its
complexification.  This induces a remarkably rich structure on many
associated invariants of $V$.  A starting point for all of this is the
fact that $\M_g$ is a complex algebraic variety.  Since, at the end of the
day, Hain's presentation of ${\cal L}(\T_g)\otimes \R$ is rather
simple, it is natural to pose the following.

\begin{problem}
\label{problem:hain1}
Give an elementary, purely combinatorial-topological and
group-theoretic, proof of Hain's theorem.
\end{problem}

It seems that a solution to Problem \ref{problem:hain1} will likely
require us to advance our understanding of $\T_g$ in new ways.  It may
also give a hint towards attacking the following problem, where mixed
Hodge theory does not apply.

\begin{problem}[Hain for $\Aut(F_n)$]
Give an explicit finite presentation for the Malcev Lie Algebra 
${\cal L}({\rm IA}_n)$, where ${\rm IA}_n$ is the group of automorphisms
of the free group $F_n$ acting trivially on $H_1(F_n,\Z)$.
\end{problem}

There is a great deal of interesting information at the prime $2$ which
Hain's theorem does not address, and indeed which remains largely
unexplored.  While Hain's theorem tells us that reduction$\mod 2$ gives
us a large subalgebra of ${\cal L}_2(\T_{g,1})$ coming from ${\cal
L}_0(\T_{g,1})$, the Lie algbera ${\cal L}_2(\T_{g,1})$ over $\F_2$ is
much bigger.  This can already be seen from (\ref{eq:h1:torelli}).  
As noted above, the
$2$-torsion $B_2$ exists for ``deeper'' reasons than the other piece of
$H_1(\T_{g,1},\Z)$, as it comes from the Rochlin invariant as opposed to
pure algebra. Indeed, for the analogous ``Torelli group'' $\IA_n$ for
$\Aut(F_n)$, the corresponding ``Johnson homomorphism'' gives all the
first cohomology.  Thus the $2$-torsion in $H^1(\T_g,\Z)$ is truly
coming from $3$-manifold theory.

\begin{problem}[Malcev mod $2$]
Give an explicit finite 
presentation for the $\F_2$-Lie algebra ${\cal L}_2(\T_{g,1})$.  
\end{problem}

We can also build a Lie algebra using the Johnson filtration.  Let 
$${\mathfrak h}_g:=\bigoplus_{k=1}^\infty\frac{\displaystyle
\T_g(k)}{\displaystyle \T_g(k+1)}\otimes \R$$ 
Then ${\mathfrak h}$ is a real Lie algebra.  In \S 14 of \cite{Ha3}, 
Hain proves  that the Johnson filtration
is not cofinal with the lower central series of $\T_g$.  He also relates
${\mathfrak h}_g$ to ${\mathfrak t}_g$.  The 
following basic question remains open.  

\begin{question}[Lie algebra for the Johnson filtration]
Is ${\mathfrak h}_g$ a finitely presented Lie algebra? If so, give an 
explicit finite presentation for it.
\end{question}

\subsection{Low-dimensional homology of principal congruence subgroups}

Recall that the {\em level $L$ congruence subgroup} $\Mod_g[L]$ is
defined to be the subgroup of $\Mod_g$ which acts trivially on
$H_1(\Sigma_g,\Z/L\Z)$.  This normal subgroup has finite index; indeed
the quotient of $\Mod_g$ by $\Mod_g[L]$ is the finite symplectic group
$\Sp(2g,\Z/L\Z)$.  When $L\geq 3$ the group $\Mod_g[L]$ is torsion free,
and so the corresponding cover of moduli space is actually a manifold.  
These manifolds arise in algebraic geometry as they 
parametrize so-called ``genus $g$ curves with level $L$ structure''; see
\cite{Ha2}.

\begin{problem}
Compute $H_1(\Mod_g[L];\Z)$.
\end{problem}

McCarthy and (independently) Hain proved that $H_1(\Mod_g[L],\Z)$ is 
finite for $g\geq 3$; see, e.g. Proposition 5.2 of
\cite{Ha2}\footnote{Actually, Hain proves a much stronger result,
computing $H^1(\Mod_g[L],V)$ for $V$ any finite-dimensional 
symplectic representation.}.  As
discussed in \S 5 of \cite{Ha2}, the following conjecture would imply
that the (orbifold) Picard group for the moduli spaces of level $L$
structures has rank one; this group is finitely-generated by the Hain
and McCarthy result just mentioned.

\begin{conjecture}[Picard number one conjecture for level $L$ structures]
\label{conjecture:pic1}
Prove that $H_2(\Mod_g[L];\Q)=\Q$ when $g\geq 3$.  
More generally, compute $H_2(\Mod_g[L];\Z)$ for all $g\geq 3, L\geq 2$.   
\end{conjecture}

Harer \cite{Har2} proved this conjecture in the case $L=1$.  This
generalization was stated (for Picard groups) as Question 7.12 in
\cite{HL}.  The case 
$L=2$ was claimed in \cite{Fo}, but there is apparently 
an error in the proof.  At this point even 
the $(g,L)=(3,2)$ case is open.  

Here is a possible approach to Conjecture \ref{conjecture:pic1} for
$g\geq 4$.  First
note that, since $\Mod_g[L]$ is a finite index subgroup of the finitely
presented group $\Mod_g$, it is finitely presented.  As we have a lot of
explicit information about the finite group $\Sp(2g,\Z/L\Z)$, it seems 
possible in principle to answer the following, which is also a test of 
our understanding of $\Mod_g[L]$.

\begin{problem}[Presentation for level $L$ structures]
\label{problem:level:presentation}
Give an explicit finite presentation for $\Mod_g[L]$.
\end{problem}

Once one has such a presentation, it seems likely that it would fit well
into the framework of Pitsch's proof \cite{Pi} that 
$\rank(H_2(\Mod_{g,1},\Z))\leq 1$ for $g\geq 4$.  Note that Pitsch's
proof was extended to
punctured and bordered case by Korkmaz-Stipsicsz; see \cite{Ko}.  
What Pitsch does is to
begin with an explicit, finite presentation of $\Mod_{g,1}$, and then to
apply Hopf's formula for groups $\Gamma$ presented as the quotient 
of a free group $F$ by the normal closure $R$ of the relators:
\begin{equation}
\label{eq:Hopf}
H_2(\Gamma,\Z)=\frac{\displaystyle R\cap [F,F]}{\displaystyle [F,R]}
\end{equation}
In other words, elements of $H_2(\Gamma,\Z)$ come precisely from
commutators which are relators, except for the trivial ones.  Amazingly,
one needs only write the form of an arbitrary element of the numerator
in (\ref{eq:Hopf}), and a few tricks reduces the computation of 
(an upper bound for) space of solutions to computing the rank of an
integer matrix.  In our case this approach seems feasible, especially
with computer computation, at least for
small $L$.  Of course one hopes to find a general pattern.

\section{Linear and nonlinear representations of $\Mod_g$}

While for $g>2$ it is not known whether or not $\Mod_g$ admits a faithful,
finite-dimensional  linear
representation, there are a number of known linear and nonlinear
representations of $\Mod_g$ which are quite useful, 
have a rich internal structure, and connect to other problems.  In this
section we pose a few problems about some of these.

\subsection{Low-dimensional linear representations}

It would be interesting to classify all irreducible complex 
representations $\psi:\Mod_g\to
\GL(m,\C)$ for $m$ sufficiently small compared to $g$.  This was done
for representations of the $n$-strand braid group for $m\leq n-1$ 
by Formanek \cite{For}.  There are a number of
special tricks using torsion in $\Mod_g$ and so, as with many questions
of this type, one really wants to understand low-dimensional irreducible
representations of the (typically torsion-free) 
finite index subgroups of $\Mod_g$.  It is proven in \cite{FLMi} that no
such {\em faithful} representations exist for $n<2\sqrt{g-1}$.  

One question is to determine if the standard representation
on homology $\Mod_g\to\Sp(2g,\Z)$ is minimal in some sense.  Lubotzky
has found finite index subgroups $\Gamma<\Mod_g$ and 
surjections $\Gamma\to \Sp(2g-2,\Z)$.  I believe that it should be
possible to prove that representations of such $\Gamma$ in low degrees 
must have traces which are algebraic integers.  This
problem, and various related statements providing constraints on
representations, reduce via now-standard methods to the problem of
understanding representations $\rho:\Mod_g\to \GL(n,K)$, where $K$ is a
discretely valued field such as the $p$-adic rationals.  The group 
$\GL(n,K)$ can be realized as a group of isometries of the 
corresponding Bruhat-Tits affine building; this is a nonpositively
curved (in the ${\rm CAT}(0)$ sense), $(n-1)$-dimensional 
simplicial complex.  The general problem then becomes:

\begin{problem}[Actions on buildings]
\label{problem:buildings}
Determine all isometric actions $\psi:\Mod_g\to\Isom(X^n)$, where $X^n$ is 
an $n$-dimensional Euclidean building, and $n$ is sufficiently small
compared to $g$.
\end{problem}

For example, one would like conditions under which $\psi$ has a
{\em global fixed point}, that is, a point $x\in X^n$ such that
$\psi(\Mod_g)\cdot x=x$.  One method to attack this problem is the 
so-called ``Helly technique'' introduced in \cite{Fa4}.  
Using standard ${\rm CAT}(0)$ methods, one can show that each Dehn
twist $T_\alpha$ in $\Mod_g$ has a nontrivial fixed set $F_\alpha$ under
the $\psi$-action; $F_\alpha$ is necessarily convex.  Considering the
nerve of the collection $\{F_\alpha\}$ gives a map 
${\cal C}_g\to X^n$ from the complex of curves to $X^n$.  Now ${\cal
C}_g$ has the homotopy type of a wedge of spheres (see, e.g.,
\cite{Iv1}), while $X^n$ is contractible.  Hence the 
spheres in the nerve must be filled in, which gives that many more elements 
$\psi(T_\alpha)$ have common fixed points.  The problem now is to
understand in an explicit way the spheres inside ${\cal C}_g$.

\subsection{Actions on the circle}

It was essentially known to Nielsen that $\Mod_{g,1}$ acts faithfully by
orientation-preserving homeomorphisms on the circle.  Here is how this
works: for $g\geq 2$ any homeomorphism $f\in\Homeo^+(\Sigma_g)$ lifts to
a quasi-isometry $\widetilde{f}$ of the hyperbolic plane $\hyp^2$.  Any
quasi-isometry of $\hyp^2$ takes geodesic rays to a uniformly bounded
distance from geodesic rays, thus inducing a map $\partial
\widetilde{f}:S^1\to S^1$ on the boundary of infinity of $\hyp^2$, 
which is easily checked to be a homeomorphism,
indeed a quasi-symmetric homeomorphism.  If $h\in\Homeo^+(\Sigma_g)$ is
homotopic to $f$, then since homotopies are compact one sees directly
that $\widetilde{h}$ is homotopic to $\widetilde{f}$, and so these maps
are a bounded distance from each other in the sup norm.  In particular
$\partial \widetilde{h}=\partial\widetilde{f}$; that is, 
$\partial\widetilde{f}$ 
depends only the homotopy class of $f$.  It is classical that
quasi-isometries are determined by their boundary values, hence  
$\partial\widetilde{f}=Id$ only when $f$ is homotopically trivial.  Now
there are $\pi_1\Sigma_g$ choices for lifting any such $f$, so the
group $\Gamma_g\subset \Homeo^+(S^1)$ 
of all lifts of all homotopy classes of 
$f\in\Homeo^+(\Sigma_g)$ gives
an exact sequence
\begin{equation}
\label{eq:group:of:lifts}
1\to\pi_1\Sigma_g\to \Gamma_g\to\Mod_{g,1}\to 1
\end{equation}
Since each element of $\Mod_{g,1}$ fixes a common marked point on
$\Sigma_g$, there is a canonical way to choose a lift of each $f$; that
is, we obtain a section $\Mod_{g,1}\to\Gamma_g$ splitting
(\ref{eq:group:of:lifts}). In particular we have an injection 
\begin{equation}
\label{eq:circle1}
\Mod_{g,1}\hookrightarrow \Homeo^+(S^1)
\end{equation}

This inclusion provides a so-called (left) circular ordering on
$\Mod_{g,1}$ - see \cite{Cal}.  
Note that no such inclusion as in (\ref{eq:circle1}) 
exists for $\Mod_g$ since any finite
subgroup of $\Homeo^+(S^1)$ must be cyclic\footnote{One can see this by 
averaging any
Riemannian metric on $S^1$ by the group action.}, 
but $\Mod_g$ has noncyclic finite subgroups.  

The action given by (\ref{eq:circle1}) gives a dynamical 
description of $\Mod_{g,1}$ via its action on $S^1$.  For example, each 
pseudo-Anosov in $\Mod_{g,1}$ acts on $S^1$ with finitely many fixed
points, with alternating sources and sinks as one moves around the
circle (see, e.g., Theorem 5.5 of \cite{CB}). There is an intrinsic
non-smoothness to this action, and indeed in \cite{FF} it is proven that any
homomorphism $\rho:\Mod_{g,1}\to \Diff^2(S^1)$ has trivial
image; what one really wants to prove is that no finite index subgroup of
$\Mod_g$ admits a faithful $C^2$ action on $S^1$.  
It would be quite interesting to prove that 
the action (\ref{eq:circle1}) is canonical, in the following sense.

\begin{question}[Rigidity of the $\Mod_{g,1}$ action on $S^1$]
Is any faithful action $\rho:\Mod_{g,1}\to\Homeo^+(S^1)$ 
conjugate in $\Homeo^+(S^1)$ to the standard action, given in
(\ref{eq:circle1})? What about the same question for finite index
subgroups of $\Mod_{g,1}$?
\end{question}

Perhaps there is a vastly stronger, topological dynamics characterization
of $\Mod_{g,1}$ inside $\Homeo^+(S^1)$, in the style of the Convergence
Groups Conjecture (theorem of Tukia, Casson-Jungreis and Gabai), with
``asymptotically source -- sink'' being replaced here by ``asymptotically
source -- sink -- $\cdots$ -- source -- sink'', or some
refinement/variation of this.

Now, the group of lifts of elements of $\Homeo^+(S^1)$ to homeomorphisms
of $\R$ gives a central extension 
$$
1\to\Z\to\widetilde{\Homeo}(S^1)\to\Homeo^+(S^1)\to 1
$$
which restricts via (\ref{eq:circle1}) to a central extension
\begin{equation}
\label{eq:central:ext}
1\to\Z\to\widetilde{\Mod_{g,1}}\to\Mod_{g,1}\to 1
\end{equation}

Note that $\Mod_{g,1}$ has torsion.  Since $\widetilde{\Homeo}(S^1)\subset
\Homeo^+(\R)$ which clearly has no torsion, it follows that
(\ref{eq:central:ext}) does not split.  In particular the extension
(\ref{eq:central:ext}) gives a nonvanishing class 
$\xi\in H^2(\Mod_{g,1},\Z)$.  Actually, it is not hard to see that
$\xi$ is simply the ``euler cocycle'', which assigns to any pointed map
$\Sigma_h\to\M_g$ the euler class of the pullback bundle of the
``universal circle bundle'' over $\M_g$.

The torsion in $\Mod_{g,1}$ and in $\Mod_g$ preclude each from having a
left-ordering, or acting faithfully on $\R$.  As far as we know
this is the only obstruction; it disappears when one passes to
appropriate finite index subgroups.

\begin{question}[orderability]
Does $\Mod_g, g\geq 2$ have some finite index subgroup which acts
faithfully by homeomorphisms on $S^1$? Does either $\Mod_g$ or
$\Mod_{g,1}$ have a finite index subgroup which acts faithfully by
homeomorphisms on $\R$? 
\end{question}

Note that Thurston proved that braid groups are orderable.  Since $\T_g$
and $\T_{g,1}$ are residually torsion-free nilpotent, they are 
isomorphic to a
subgroup of $\Homeo^+(\R)$; in fact one can show that
(\ref{eq:central:ext}) splits when restricted to $\T_{g,1}$.  
On the other hand, Witte \cite{Wi} 
proved that no finite index subgroup
of $\Sp(2g,\Z)$ acts faithfully by homoeomorphisms 
on $S^1$ or on $\R$.  

\medskip
\noindent
{\bf Non-residual finiteness of the universal central extension. }
The Lie group $\Sp(2g,\R)$ has infinite cyclic fundamental group.  
Its universal cover $\widetilde{\Sp(2g,\R)}$ gives a central
extension
$$1\to\Z\to\widetilde{\Sp(2g,\R)}\to\Sp(2g,\R)\to 1$$
which restricts to a central extension 
\begin{equation}
\label{eq:deligne}
1\to\Z\to\widetilde{\Sp(2g,\Z)}\to\Sp(2g,\Z)\to 1
\end{equation}
The cocycle $\zeta\in H^2(\Sp(2g,\Z),\Z)$ defining the extension
(\ref{eq:deligne}) is nontrivial and bounded; 
this comes from the fact that it is
proportional to the Kahler class on the corresponding locally symmetric
quotient (which is a $K(\pi,1)$ space).  Deligne proved in \cite{De} 
that $\widetilde{\Sp(2g,\Z)}$ is {\em not
residually finite}.  Since there is an obvious surjection of exact
sequences from (\ref{eq:circle1}) to (\ref{eq:deligne}), and since both
central extensions give a bounded cocycle, one begins to
wonder about the following.

\begin{question}[(Non)residual finiteness]
Is the (universal) central extension $\widetilde{\Mod_{g,1}}$ of
$\Mod_{g,1}$ residually finite, or not?
\end{question}

Note that an old result of Grossman states that 
$\Mod_g$ and $\Mod_{g,1}$ are both residually finite.  The group 
$\Sp(2g,\Z)$ is easily seen to be residually finite; indeed the
intersection of all congruence subgroups of $\Sp(2g,\Z)$ is trivial.

\subsection{The sections problem}

Consider the natural projection $\pi:\Homeo^+(\Sigma_g)\to \Mod_g$, and
let $H$ be a subgroup of $\Mod_g$.  We say that $\pi$ {\em has a section over
$H$} if there exists a homomorphism $\sigma:\Mod_g\to
\Homeo^+(\Sigma_g)$ so that $\pi\circ \sigma =\Id$.  This means
precisely that $H$
has a section precisely when it can be realized as a group of
homeomorphisms, not just a group of homotopy classes of
homeomorphisms. The general problem is then the following.

\begin{problem}[The sections problem]
\label{problem:sections}
Determine those subgroups $H\leq \Mod_g$ for which $\pi$ has a section
over $H$.  Do this as well with $\Homeo^+(\Sigma_g)$ replaced by 
various subgroups, such as $\Diff^r(S)$ with 
$r=1,2,\ldots ,\infty, \omega$; similarly for the group of area-preserving
diffeomorphisms, quasiconformal homeomorphisms, etc.. 
\end{problem}

Answers to Problem \ref{problem:sections} are known in a number of
cases.

\begin{enumerate}

\item When $H$ is free then sections clearly always exist over $H$. 

\item Sections to $\pi$ exist over free abelian $H$, even when
restricted to $\Diff^\infty(\Sigma_g)$.  This
is not hard to prove, given the classification by
Birman-Lubotzky-McCarthy \cite{BLM} of abelian subgroups of $\Mod_g$.

\item Sections exist over any finite group $H<\Mod_g$, even when
restricted to $\Diff^\omega(\Sigma_g)$.  
This follows from the Nielsen Realization Conjecture, proved
 by Kerckhoff \cite{Ke}, which states that any such $H$ acts as a group
of automorphisms of some genus $g$ Riemann surface.

\item In contrast, Morita showed (see, e.g., \cite{Mo5}) 
that $\pi$ does not have a section with image in $\Diff^2(\Sigma_2)$ 
over all of $\Mod_g$ when $g\geq 5$.  The $C^2$ assumption is used in a
crucial way since Morita uses a putative section to build a codimension
$2$ foliation on the universal curve over $\M_g$, to whose normal bundle
he applies the Bott vanishing theorem, contradicting nonvanishing of
a certain (nontrivial!) Miller-Morita-Mumford class.  
It seems like Morita's proof can 
be extended to finite index subgroups of $\Mod_g$.  

\item 
Markovic \cite{Mar} has recently proven that $H=\Mod_g$ does not even
have a section into $\Homeo(\Sigma_g)$, answering a well-known 
question of Thurston.  
\end{enumerate}

As is usual when one studies representations of a discrete group
$\Gamma$, one really desires a theorem about all finite 
index subgroups of $\Gamma$.  One reason for this is that 
the existence of torsion and special relations in a group $\Gamma$ 
often highly constrains its possible representations.  
Markovic's proof in \cite{Mar} uses both torsion and the braid relations
in what seems to be an essential way; these both disappear in most
finite index subgroups of $\Mod_g$.  
Thus it seems that a new idea is needed to answer the
following.

\begin{question}[Sections over finite index subgroups]
Does the natural map $\Homeo^+(\Sigma_g)\to\Mod_g$ have a section over 
a finite index subgroup of $\Mod_g$, or not?
\end{question}

Of course the ideas in \cite{Mar} are likely to be pertinent.  
Answers to Problem \ref{problem:sections} even for specific subgroups 
(e.g. for $\T_g$ or more generally $\T_g(k)$) would be interesting.  
It also seems reasonable to believe that the existence of sections is
affected greatly by the degree of smoothness one requires.   

Instead of asking for {\em sections} in the above questions, one can ask more
generally whether there are any actions of $\Mod_g$ on $\Sigma_g$.  

\begin{question}
Does $\Mod_g$ or any of its finite index subgroups have any faithful
action by homeomorphisms on $\Sigma_g$?
\end{question}

\section{Pseudo-Anosov theory}

Many of the problems in this section come out of joint work
with Chris Leininger and Dan Margalit, 
especially that in the paper \cite{FLM}.

\subsection{Small dilatations}

Every pseudo-Anosov mapping class $f \in \Mod_g$ has a {\em dilatation} 
$\lambda(f)\in\R$.  This number is an algebraic integer which 
records the exponential
growth rate of lengths of curves under iteration of $f$, in any fixed
metric on $S$; see \cite{Th}.  
The number $\log(\lambda(f))$ equals the minimal topological entropy of
any element in the homotopy class $f$; this minimum is realized by a
pseudo-Anosov homeomorphism representing $f$ (see \cite[Expos\'e
10]{FLP}). $\log(\lambda(f))$ is also the
translation length of $f$ as an isometry of the {\em Teichm\"uller space
of $S$} equipped with the Teichm\"uller metric.  Penner considered the 
set $$\spec(\Mod_g)=\{\log(\lambda(f)): \mbox{$f\in \Mod_g$ is
pseudo-Anosov}\} \subset \R$$ 
This set can be 
thought of as the {\em length spectrum} of $\M_g$.  
We can also consider, for various subgroups
$H<\Mod_g$, the subset $\spec(H) 
\subset
\spec(\Mod_g)$ obtained by restricting to pseudo-Anosov elements of $H$.
Arnoux--Yoccoz \cite{AY} and Ivanov
\cite{Iv2} proved that $\spec(\Mod_g)$ is discrete as a subset of $\R$.
It follows that for any subgroup $H<\Mod_g$,
the set $\spec(H)$ has a least element $L(H)$.  Penner proved in
\cite{Pe} that
$$\frac{\displaystyle\log 2}{\displaystyle 12g-12}\leq L(\Mod_g)\leq 
\frac{\displaystyle\log 11}{\displaystyle g}$$

In particular, as one increases the genus, there are pseudo-Anosovs with
stretch factors arbitrarily close to one.  In contrast to our
understanding of the asymptotics of $L(\Mod_g)$, we still do not know
the answer to the following question, posed by McMullen.

\begin{question}
\label{question:curt3}
Does $\lim_{g\to\infty}gL(\Mod_g)$ exist?
\end{question}

Another basic open question is the following.

\begin{question}
Is the sequence $\{L(\Mod_g)\}$ monotone decreasing? strictly so?
\end{question}

Explicit values of $L(\Mod_g)$ are known only when $g=1$.  In this case
one is simply asking for the minimum value of the largest root of a
polynomial as one varies over all integral 
polynomials $x^2-bx+1$ with $b\geq
3$.  This is easily seen to occur when $b=3$. For $g=2$ Zhirov \cite{Zh}
found the smallest dilatation for pseudo-Anosovs with orientable
foliation.  It is not clear if this should equal $L(\Mod_2)$ or not.

\begin{problem}
\label{problem:smallest:dil}
Compute $L(\Mod_g)$ explicitly for small $g\geq 2$.
\end{problem}

In principle $L(\Mod_g)$ can be computed for any given $g$.  The point is
that one can first bound the degree of $L(\Mod_g)$, then give bounds on 
the smallest possible value $\lambda(\alpha)$, where $\alpha$ ranges
over all algebraic integers of a fixed range of degrees, and
$\lambda(\alpha)$ denotes the largest root of the minimal polynomial of
$\alpha$.  One then finds all train tracks on $\Sigma_g$, and starts to
list out all pseudo-Anosovs.  It is possible to give bounds for when
the dilatations of these become large.  Now one tries to match up the
two lists just created, to find the minimal dilatation pseudo-Anosov on
$\Sigma_g$.  Of course actually following out this procedure, even for
small $g$, seems to be impracticable.

\begin{question}
Is there a unique (up to conjugacy) minimal dilation pseudo-Anosov 
in $\Mod_g$?
\end{question}

Note that this is true for $g=1$; the unique minimum is realized
by the conjugacy class of the matrix $\left(\begin{array}{cc}
2&1\\
1&1
\end{array}\right)$.

Here is a natural refinement of the problem of finding $L(\Mod_g)$.  
Fix a genus $g$.  Fix a possible $r$-tuple $(k_1,\ldots ,k_r)$ of
{\em singularity data} for $\Sigma_g$.  By this we mean to consider
possible foliations with $r$ singularities with $k_1,\ldots ,k_r$
prongs, respectively.  For a fixed $g$, there are only finitely many
possible tuples, as governed by the Poincare-Hopf index theorem.  
Masur-Smillie \cite{MS} proved that, for every admissible tuple, there is some
pseudo-Anosov on $\Sigma_g$ with stable foliation having the given
singularity data.  Hence the following makes sense.

\begin{problem}[Shortest Teichm\"{u}ller loop in a stratum]
\label{problem:strata:systoles}
For each fixed $g\geq 2$, and for each $r$-tuple as above, give upper
and lower bounds for 
$$\lambda_g(k_1,\ldots ,k_r):=\inf\{\log\lambda(f): \mbox{$f\in \Mod_g$
whose stable foliation has data $(k_1,\ldots ,k_r)$}\}$$
\end{problem}

This problem is asking for bounds for 
the shortest Teichm\"{u}ller loop lying in a 
given substratum in moduli space (i.e. the projection in $\M_g$
of the corresponding substratum in the cotangent bundle).
$$L(\Mod_g)=\min\{\lambda_g(k_1,\ldots ,k_r)\}$$
where the $\min$ is taken of all possible singularity data.  

\subsection{Multiplicities}

The set $\spec(\Mod_g)$ has unbounded multiplicity; that is, given any
$N>0$, there exists $r\in\spec(\Mod_g)$ such that there are at least $n$
conjugacy classes $f_1,\ldots f_n$ of pseudo-Anosovs in $\Mod_g$ having 
$\log(\lambda(f_i))=r$.  The reason for this is that $\M_g$ contains
isometrically embedded finite volume hyperbolic $2$-manifolds $X$, e.g. the
so-called {\em Veech curves}, and any such $X$ has (hyperbolic) length
spectrum of unbounded mulitplicity.  

A related mechanism which produces
length spectra with unbounded multiplicities is the {\em Thurston 
representation}.  This gives, for a pair of curves $a,b$ on $\Sigma_g$
whose union fills $\Sigma_g$, an injective representation
$\rho:<T_a,T_b> \to \PSL(2,\R)$ with the following properties:
$\image(\rho)$ is discrete; each element of $\image(\rho)$ is either 
pseudo-Anosov or is a power of $\rho(T_a)$ or $\rho(T_b)$; and 
$\spec(<T_a,T_b>)$ is essentially the length spectrum of the quotient of
$\hyp^2$ by $\image(\rho)$. Again it follows that $\spec(<T_a,T_b>)$ has
unbounded multiplicity.  Since one can find $a,b$ as above, each of
which is in addition separating, 
it follows that $\spec(\T_g)$ and even $\spec(\T_g(2))$ have
unbounded multiplicity.  

\begin{question}
Does $\spec(\T_g(k))$ have bounded multiplicity for $k\geq 3$?
\end{question}

One way to get around unbounded multiplicities is to look
at the {\em simple} length spectrum, which is the subset of
$\spec(\Mod_g)$ coming from pseudo-Anosovs represented by 
simple (i.e. non-self-intersecting) geodesic loops in $\M_g$.

\begin{question}[Simple length spectrum]
\label{question:multiplicity1}
Does the simple length spectrum of $\M_g$, endowed with 
the Teichm\"{u}ller metric, have bounded multiplicity? If so, how does
the bound depend on $g$?
\end{question}

Of course this question contains the corresponding question for (many)
hyperbolic surfaces, which itself is still open. These questions also
inspire the following.

\begin{problem}
Give an algorithm which tells whether or not any given pseudo-Anosov is
represented by a simple closed Teichm\"{u}ller geodesic, and also
whether or not this geodesic lies on a Veech curve.
\end{problem}
                              
Note that the analogue of Question \ref{question:multiplicity1} is not
known for hyperbolic surfaces, although it is true for a generic set of
surfaces in $\M_g$.

\subsection{Special subgroups}
In \cite{FLM} we provide evidence for the
principle that algebraic complexity implies dynamical complexity.  A 
paradigm for this is the following theorem.

\begin{theorem}[\cite{FLM}]
\label{theorem:torelli} For $g \geq 2$, we have
$$ .197 < L(\T_g) < 4.127$$
\end{theorem}

The point is that $L(\T_g)$ is universally bounded above and below, 
independently of $g$.  We extend this kind of universality to every term
of the Johnson filtration, as follows.

\begin{theorem}[\cite{FLM}]
\label{theorem:johnson filtration}
Given $k \geq 1$, there exist $M(k)$ and $m(k)$, where $m(k) \to
\infty$ as $k \to \infty$, so that
\[ m(k) < L(\T_g(k)) < M(k) \]
for every $g \geq 2$.
\end{theorem}

\begin{question}
Give upper and lower bounds for $L(\T_g(k))$ for all $k\geq 2$ which are
of the same order of magnitude.  
\end{question}

In \cite{FLM} bounds on $L(H)$ are given for various special classes of
subgroups $H<\Mod_g$.  It seems like there is much more to explore in
this direction.  One can also combine these types of questions with
problems such as Problem \ref{problem:strata:systoles}.

\subsection{Densities in the set of dilatations} 

Let P be a property which pseudo-Anosov homeomorphisms might or might
not have. For example, P might be the property of lying in a fixed
subgroup of $H<\Mod_g$, such as $H=\T_g$; or P might be the property of
having dilatation which is an algebraic integer of a fixed degree.  
One can then ask the natural question: how commonly do the dilatations
of elements with P arise in $\spec(\Mod_g)$?  

To formalize this, recall that the {\em (upper) density} $d^\ast(A)$ of
a subset $A$ of the natural numbers $\N$ is defined as
$$d^\ast(A):=\limsup_{N\to\infty}\frac{\displaystyle \# A\cap
[0,n]}{\displaystyle n}$$ This notion can clearly be extended from $\N$
to any countable ordered set ${\cal S}$ once an order-preserving
bijection ${\cal S}\to \N$ is fixed.

Now fix $g\geq 2$, and recall that $\spec(\Mod_g)\subset \R^+$ is defined
to be the set of ($\log$s of) dilatations of pseudo-Anosov
homeomorphisms of $\Mod_g$.  The set $\spec(\Mod_g)$ comes with 
a natural order $\lambda_1<\lambda_2<\cdots $, which determines a fixed
bijection $\spec(\Mod_g)\to\N$.  If we wish to keep track of all
pseudo-Anosovs, and not just their dilatations, we can simply consider
the (total) ordering on the set of all pseudo-Anosovs ${\cal P}_g$ given
by $f\leq g$ if $\lambda(f)\leq \lambda(g)$.  

\begin{question} 
For various subgroups $H<\Mod_g$, compute the density of $\spec(H)$ in 
$\spec(\Mod_g)$ and the density of $H\cap {\cal P}_g$ in ${\cal P}_g$.  
In particular, what is the density of
$\spec(\Mod_g[L])$ in $\spec(\Mod_g)$? What about $H=\T_g(k)$ with
$k\geq 1$?
\end{question}

\noindent
Dept. of Mathematics, University of Chicago\\
5734 University Ave.\\
Chicago, Il 60637\\
E-mail: farb@math.uchicago.edu

\end{document}